\documentclass[11pt]{amsart}
\def \ba {\begin {eqnarray*} }
\def \ea {\end {eqnarray*} }
\def \beq {\begin {eqnarray}}    
\def \eeq {\end {eqnarray}}         
\def \supp {\hbox{supp}\,}

\def \e {\varepsilon}
\def \p {\partial}
\def \a {\alpha}
\def \la{{\lambda}}

\def\R{\mathbb R}

\def\C{\mathbb C}
\def\Z{\mathbb Z}
\def\O{{\mathcal O}}
\def\M{{\mathcal M}}

\newcommand{\nntext}{}
\newcommand{\ntekst}{}

\newtheorem{theorem}{Theorem}[section]
\newtheorem{lemma}[theorem]{Lemma}           
\newtheorem{cor}[theorem]{Corollary}

\theoremstyle{definition}
\newtheorem{definition}[theorem]{Definition}

\newtheorem{condition}[theorem]{Condition}
\theoremstyle{remark}
\newtheorem{remark}[theorem]{Remark}

\numberwithin{equation}{section}

\begin{document}
\baselineskip 13pt

\author{Hiroshi Isozaki}
\address{Hiroshi Isozaki, Institute of Mathematics \\
University of Tsukuba,
Tsukuba, 305-8571, Japan}
\author{Yaroslav Kurylev}
\address{Yaroslav Kurylev, Department of Mathematics \\
University College of London, United Kingdom}
\author{Matti Lassas}
\address{Matti Lassas, Department of Mathematics and Statistics\\
University of Helsinki, Finland}

\title[Inverse scattering on
asymptotically hyperbolic surfaces]{Conic singularities, generalized scattering matrix, and inverse scattering on
asymptotically hyperbolic surfaces}
\date{August 3, 2011}

\maketitle

\begin{abstract}
We consider an inverse problem associated with some 2-dimensional non-compact surfaces with conical singularities, cusps
{\ntekst and regular ends.}  Our 
motivating example is a Riemann surface $\mathcal M = \Gamma\backslash{\bf H}^2$ associated with a Fuchsian group of the 1st kind $\Gamma$ 
containing parabolic elements.  $\mathcal M$ is then non-compact, and has a finite number of cusps and elliptic singular points, which is 
regarded as a hyperbolic orbifold. We introduce a class of Riemannian surfaces with conical singularities on its finite part, having cusps 
{\ntekst and regular ends} at infinity, 
whose metric is asymptotically hyperbolic. By observing solutions of the Helmholtz equation at the cusp, we define a generalized S-matrix. We then 
show that this generalized S-matrix determines the Riemannian metric and the  structure of conical singularities.
\end{abstract}


\section{Introduction}


\subsection{Assumptions on the manifold}
Throughout this paper $S^{r}$ denotes the circle of radius $r$, which is identified 
with ${\bf R}/(2\pi r{\bf Z})$: 
\ba
S^{r} = \{(x_1,x_2) \, &;& x_1 + ix_2 = re^{(ix/r)}=
(r \cos(\theta),\, r\sin(\theta)),\\
& &0 \leq x \leq 2\pi r, \,\, 0 \leq \theta \leq 2\pi\},
\nonumber
\ea
with an obvious identification of $\theta=0$ and $\theta=2\pi$. Thus, to write a function
on $S^r$ we would write it as $f(x),\, x \in S^r$, or $f(\theta),\, \theta \in [0, 2\pi]$
or $\theta \in \R$, assuming $2\pi-$periodicity.

 We consider a $2$-dimensional orientable connected $C^{\infty}$-surface 
 without boundary
 ${\mathcal M}$, which is 
 written as a union of open sets:
\begin{equation}
{\mathcal M} = {\mathcal K}\cup{\mathcal M}_1\cup\cdots\cup{\mathcal M}_N
\label{1.0}
\end{equation}
satisfying the following assumptions. 

\medskip
\noindent
{\bf (A-1)} {\it ${\overline{\mathcal K}}$ {\it is compact}, and there exists a finite set $\mathcal M_{sing} \subset \mathcal K$ such that 
 $\mathcal M\setminus \mathcal M_{sing}$
 is equipped with a $C^{\infty}$-Riemannian metric $g$}. 

\medskip
\noindent
{\bf (A-2)} {\it For any $p \in \mathcal M_{sing}$, there exist a constant $\epsilon_p > 0$ and local coordinates $(r,\theta) \in (0,\epsilon_p)\times [0, 2\pi]$ around $p$ such that $r = 0$ corresponds to $p$ and the Riemannian metric $g$ takes the form}
\begin{equation}
(ds)^2 = (dr)^2 + C_pr^2(1 + h_p(r,\theta))(d\theta)^2,
\label{S1metricaroundsingular}
\end{equation}
where

\medskip
\noindent
(A-2-1) $C_p$ {\it is a positive constant such that} $C_p \neq 1$, 

\medskip
\noindent
(A-2-2) $h_p(r,\theta) \in C^{\infty}((0,\epsilon_p)\times [0, 2\pi])$,

\medskip
\noindent
(A-2-3) {\it As $r \to 0$, $h_p(r,\theta) \to 0$ uniformly with respect to 
$\theta \in [0, 2\pi]$.}

\medskip
\noindent
{\bf (A-3)} $\ $ {\it There exists $\mu \geq 1$ such that for 
$1 \leq i \leq \mu$, $\mathcal M_i$ is isometric to 
$S^{r_i}\times (1,\infty),\, r_i >0,$ equipped with the metric}
\begin{equation}
ds^2 = \frac{(dx)^2 + (dy)^2}{y^2}.
\nonumber
\end{equation}

\medskip
\noindent
{\bf (A-4)} $\ $ {\it For $\mu+1 \leq i \leq N$, ${\mathcal M}_i$ is diffeomorphic to $S^{r_i}\times(0,1)$, $r_i > 0$, and the metric on $\mathcal M_i$ has the following form :
\begin{equation}
ds^2 = y^{-2}\left((dy)^2 + (dx)^2 + A(x,y,dx,dy)\right),
\nonumber
\end{equation}
\begin{equation}
A(x,y,dx,dy) = a_1(x,y)(dx)^2 + a_2(x,y)dxdy + a_3(x,y)(dy)^2,
\nonumber
\end{equation}
where  $a_i(x,y)$ $(i = 1,2,3)$ satisfies the following condition
\begin{equation}
|\partial_x^{\alpha}\big(y\partial_y\big)^n\, a_i(x,y)| \leq C_{\alpha n}(1 + |\log y|)^{-n-1-\epsilon_0}, \quad \forall \alpha, n,
\nonumber
\end{equation}
for some $\epsilon_0 > 0$.}

\medskip
We say that under the above assumption (A-2), the metric $g$ has a {\it conical singularity} at $p \in \mathcal M_{sing}$. The part ${\mathcal M}_i$, $1 \leq i \leq \mu$, will be called a {\it cusp}. (This is a little abuse of the standard terminology).
Since $\mu \geq 1$, $\mathcal M$ has at least one cusp. If $\mu = N$, all the ends have a cusp. We call $\mathcal M_i$, $\mu+1 \leq i \leq N$, {\it regular part}. The metric on regular parts  are allowed to be different from each other. 

In \cite{IsKu10}, spectral theory for asymptotically hyperbolic manifolds without conical singularities is discussed, and the arguments there can be extended to the above situation. Let $\Delta_g$ be the Laplace-Beltami operator for the metric $g$, and $H$ the Friedrichs extension of $- \Delta_g - 1/4$ associated with the quadratic form $A_g[u,v]=(\nabla u,\nabla v)  - \frac{1}{4}(u,v)$
with $u,v\in {\mathcal D}(A_g)=H^1(\mathcal M)$. 
 It has continuous spectrum $\sigma_c(H) = [0,\infty)$, and the discrete spectrum $\sigma_d(H) \subset (-\infty,0)$. If at least one of the ends is regular, there is no eigenvalues in $(0,\infty)$. If all the ends are cusps, $H$ may have embedded eigenvalues in $(0,\infty)$, which are discrete with possible accumulation points $0$ and $\infty$. 


\subsection{Inverse scattering from regular ends}
An important notion to describe the spectral properties of $H$ is the S-matrix. Usually, it 
is introduced by observing the asymptotic behavior, as time tends to $\pm \infty$,  of 
solutions to the time-dependent Schr{\"o}dinger equation or the wave equation on 
$\mathcal M$, i.e. $S = W_+^{\ast}W_-$, 
where $W_{\pm} = {\rm s-lim}_{t\to\pm\infty}e^{itH}e^{-itH_0}$, or 
$W_{\pm} = {\rm s-lim}_{t\to \pm\infty}e^{it\sqrt{H}}e^{-it\sqrt{H_0}}$, where $H_0$ is 
the {\it unperturbed} operator, to which $H$ is asymptotic at infinity. An equivalent way 
is to observe asymptotic expansions at infinity of {\it physical} solutions to the  
Helmholtz equation on $\mathcal M$. In the case of our manifold $\mathcal M$, by the 
physical solution $u$, we roughly mean that $u$ behaves like $O(y^{1/2})$ on each end. 
The ({\it physical}) S-matrix $\widehat S(k)$, $k$ being the square root of the energy of the system, is an operator valued $N\times N$ matrix, $\widehat S(k) = \big(\widehat S_{ij}(k)\big)$, where $\widehat S_{ij}(k)$ corresponds to the wave coming in from the end $\mathcal M_j$ and going out of the end $\mathcal M_i$
{(see e.g. \cite{
Bor07,BJP05,Faddeev,
Gu,
GuZw97,
GuZw95,
JSb1,
JSb2,
MM87,
Perry1,
Perry2,
Vasy}
 for variuos related results on the spectral and scattering theory for 
hyperbolic and 
asymptocally hyperbolic spaces.) Having $S-$matrix}, 
one can then talk about the inverse problem. Let us consider the case without singular points. Suppose we are given two such manifolds $\mathcal M^{(1)}$, $\mathcal M^{(2)}$, and assume  $\mathcal M_1^{(i)}$ is a regular end for $i = 1,2$.
 We also assume that for $\mathcal M^{(1)}$ and $\mathcal M^{(2)}$, the $(1,1)$ component of the associated S-matrix coincide, i.e. $\widehat S_{11}^{(1)}(k) = \widehat S_{11}^{(2)}(k)$ for all $k > 0$. If, furthermore, two ends $\mathcal M_1^{(1)}$ and $\mathcal M_1^{(2)}$ are isometric
 for large $y$, these two manifolds $\mathcal M^{(1)}$ and $\mathcal M^{(2)}$ are shown to be globally isometric (see \cite{IsKu10}). Let us note that,
 when all the ends are regular, Sa Barreto \cite{SaBa05}, see also  \cite{GuSa},  proved that,
 in the framework of scattering 
metric due to Melrose, two such manifolds are isometric, if the whole scattering matrix for all energies coincide, without assuming that one end is known to be isometric. 
The related inverse boundary value problems for compact
Riemannian manifolds can
presently be solved with fixed frequency data in
the zero energy case \cite{LU01,HM}, 
when the metric is real analytic \cite{LeU89,LTU03}, 
or when the tensor is known to be of appropriate
type up to a conformal factor \cite{DKSU09,GuT1,GuT2}. On review on the positive results
and counterexamples for these problems, see \cite{GKLU09}. 
For the resonance problem, another view point for inverse scattering, see e.g. \cite{GuZw97} and \cite{Bor07}, \cite{BP11}.


\subsection{Main result}
The problem we address here is the case in which we observe the waves coming in and 
going out of a cusp. 
Recall that the end $\mathcal M_1$ has a cusp at infinity.
Since the continuous spectrum due to the cusp is 1-dimensional, the associated S-matrix component 
$\widehat S_{11}(k)$ is a complex number, and it does not have enough information to determine 
the whole manifold. Therefore, we generalize the notion of the S-matrix. This generalized S-matrix 
was introduced in \cite{Is04} in the inverse scattering from a fixed energy  for Schr{\"o}dinger 
operators on asymptotically hyperbolic manifolds.

The Helmholtz equation has the following form in the cusp $\mathcal M_1$
\begin{equation}
- y^2(\partial_y^2 + \partial_x^2) u - \frac{1}{4}u = k^2 u,
\nonumber
\end{equation}
where $k > 0$. Passing to the Fourier series, we see that all  solutions of this equation have the asymptotic expansion 
\begin{equation}
\begin{split}
u(x,y) & \simeq a_0\,y^{\frac{1}{2}-ik} + 
\sum_{n\neq0}a_n\Big(\frac{r_1}{2\pi|n|}\Big)^{1/2}\,e^{ inx/r_1}e^{|n|y/r_1} \\
& + {b_0}\,y^{\frac{1}{2}+ik} + \sum_{n\neq0}b_n\Big(\frac{\pi r_1}{2|n|}\Big)^{1/2}\,e^{  inx/r_1}e^{-|n|y/r_1}
\end{split}
\nonumber
\end{equation}
as $y \to \infty$. We  call the operator
\begin{equation}
{\bf S}_{11}(k) : \{a_n\} \to \{b_n\}
\nonumber
\end{equation}
{\it the generalized S-matrix}, actually its (1,1) component  (see \S4 for the precise 
definition). We shall show that this generalized S-matrix determines the whole 
manifold $\mathcal M$. Namely,
suppose we are given two manifolds $\mathcal M^{(1)}, \mathcal M^{(2)}$ satisfying the assumptions 
(A-1) $\sim$ (A-4). Let $\mathcal M_{sing}^{(i)} = \{p_1^{(i)},\cdots,p_{k_i}^{(i)}\}$ be the set of singular points.

Our main result is the following.


\begin{theorem}\label{thm: main result}
Suppose we are given two manifolds $\mathcal M^{(1)}$ and $\mathcal M^{(2)}$ satifying the assumptions {\bf (A-1)} $\sim$ {\bf (A-4)}. Let the (1,1) component of the generalized scattering matrix coincide : 
$$
{\bf S}_{11}^{(1)}(k) = {\bf S}_{11}^{(2)}(k), \quad \forall k > 0, \quad  k^2\not\in \sigma_p(H^{(1)})\cup\sigma_p(H^{(2)}),
$$ 
and $r_1^{(1)}=r_1^{(2)}$.
Then there is an isometry between $\mathcal M^{(1)}$ and $\mathcal M^{(2)}$ in the following sense.

\noindent
(1) There is a homeomorphism $\Phi : \mathcal M^{(1)} \to \mathcal M^{(2)}$. \\
\noindent
(2) $\ \Phi(\mathcal M_{sing}^{(1)}) = \mathcal M_{sing}^{(2)}$.\\
\noindent
(3) $\ \Phi : \mathcal M^{(1)}\setminus\mathcal M_{sing}^{(1)} \to 
\mathcal M^{(2)}\setminus\mathcal M_{sing}^{(2)}$ is a Riemannian isometry. \\
\noindent
(4) If $p \in \mathcal M_{sing}^{(1)}$, then $C_p^{(1)}=C_{\Phi(p)}^{(2)}$ and
there is $\beta$ such that, in coordinates (A-2), we have
$h_p^{(1)}(r, \theta)=h_{\Phi(p)}^{(2)}(r, \widehat{\theta+\beta})$.

{\ntekst Here, for any $\theta \in \R$, ${\widehat \theta} \in [0, 2\pi)$ satisfies
$\theta -{\widehat \theta} \in 2 \pi \Z$.}
\end{theorem}


As will be seen from the arguments in \S 2 and \S 3, we can introduce the {\it physical} S-matrix for manifolds satisfying (A-1) $\sim$ (A-4) of this paper, and generalize the results on the inverse scattering from regular ends in \cite{IsKu10} to our case. Moreover we can also prove the same result for the inverse scattering with respect to the generalized S-matrix. Therefore, Theorem 1.1 combined with \cite{IsKu10}, implies the following theorem.


\begin{theorem}
Suppose we are given two manifolds $\mathcal M^{(1)}$ and $\mathcal M^{(2)}$ satisfying the assumptions (A-1) $\sim$ (A-4). Suppose there exist $\nu_1$ and $\nu_2$ such that the $(\nu_1,\nu_1)$ and $(\nu_2,\nu_2)$ components of the generalized S-matrices coincide:
$$
{\bf S}^{(1)}_{\nu_1\nu_1}(k) = {\bf S}^{(2)}_{\nu_2\nu_2}(k), \quad \forall k > 0, \quad k^2 \not\in \sigma_p(H^{(1)})\cup\sigma_p(H^{(2)}).
$$
Assume, furthermore, that their ends $\mathcal M^{(1)}_{\nu_1}$ and $\mathcal M^{(2)}_{\nu_2}$ are isometric. Then we have the same conclusion as in Theorem 1.1.
\end{theorem}

A good example of a surface with conical singularities is a 2-dimensional 
Riemannian orbifold, and classical examples are given by hyperbolic 
orbifolds with finite elliptic singular 
points. For example, consider $\mathcal M = \Gamma\backslash{\bf H}^2$, where $\Gamma$ is a Fuchsian group. 
As will be explained in \S2, if $\Gamma$ is a geometrically finite Fuchsian group, $\Gamma\backslash{\bf H}^2$ satisfies the assumptions (A-1) $\sim$ (A-4). Therefore, the following theorem holds.


\begin{theorem} 
Given two geometrically finite hyperbolic orbifolds $\Gamma_1\backslash{\bf H}^2$ and $\Gamma_2\backslash{\bf H}^2$, suppose there exist $\nu_1$ and $\nu_2$ such that the $(\nu_1,\nu_1)$ and $(\nu_2,\nu_2)$ components of the generalized S-matrices coincide:
$$
{\bf S}^{(1)}_{\nu_1\nu_1}(k) = {\bf S}^{(2)}_{\nu_2\nu_2}(k), \quad \forall k > 0, \quad k^2 \not\in \sigma_p(H^{(1)})\cup\sigma_p(H^{(2)}).
$$
Assume, furthermore, that their ends $\mathcal M^{(1)}_{\nu_1}$ and $\mathcal M^{(2)}_{\nu_2}$ are isometric. Then we have the same conclusion as in Theorem 1.1. Moreover, $\Phi : \mathcal M^{(1)} \to \mathcal M^{(2)}$ is an anlytic diffeomorphism, and is lifted to an orbifold isomorphism between $\mathcal M^{(1)}$ and $\mathcal M^{(2)}$.
\end{theorem}

For the notions of and geometrically finite hyperbolic orbifolds
and orbifold isomorphism, see Subsections 2.1 and 2.3.

To prove Theorem 1.1, we need to study it from two sides : the forward problem and the inverse problem. In both issues, the arguments are centered around asymptotically hyperbolic ends and singularities in the finite parts. The main ingredient of the forward problem is the spectral and scattering theory for Laplace-Beltrami operators on asymptotically hyperbolic manifolds, which two of the authors have studied in \cite{IsKu10}. Since this part does not depend on the space dimension, we shall state only the results in this paper, leaving the detailed explanations in our paper \cite{IKL11}, 
where we extend the above theorem to the higher dimensional case. 
Relations to the collapse theory of Riemannian manifolds will be discussed in \cite{KLY09}.

The crucial idea for the inverse problem part is the boundary control method. Just like our previous paper for the inverse scattering on manifolds with cylindrical ends \cite{IKL09}, we reduce the issue to the inverse boundary value problem from an artificial boundary in the end $\mathcal M_1$. 
The new ingredient in 
this paper is the argument around conic singularities based on the explicit form of  the metric (\ref{S1metricaroundsingular}).

We use a variety of notions from algebra, geometry and analysis in this paper:
Fuchsian groups, orbifolds, conical  singularities, spectral theory for self-adjoint 
operators with continuous spectrum, boundary control method. They are not 
complicated in 
themselves, however, we shall try to make the paper as readable as possible, by giving 
detailed explanations for elementary parts, sometimes referring to other papers for  
precise proofs. 
In \S 2, we recall basic facts on the Fuchsian groups, 2-dimensional hyperbolic 
orbifolds 
to explain our motivating example, and introduce the manifold with conical 
singularities. 
In \S 3, we study spectral properties of the Laplace-Beltrami operator of our manifold. 
The generalized S-matrix is defined in \S 4. We shall prove Theorem \ref{thm: main result} 
in \S 5, and Theorem 1.3 in \S 6.

The notations used in this paper are standard. For Banach spaces $\mathcal X$ and 
$\mathcal Y$, ${\bf B}(\mathcal X;\mathcal Y)$ denotes the set of all bounded linear 
operators from $\mathcal X$ to $\mathcal Y$. For a self-adjoint operator $A$ in a Hilbert 
space $\mathcal H$, $\sigma(A)$, $\sigma_p(A)$, $\sigma_{c}(A)$, $\sigma_d(A)$, 
$\sigma_{e}(A)$, $\sigma_{ac}(A)$ denote its spectrum, point spectrum (the set of all 
eigenvalues of $A$), continuous spectrum, discrete spectrum, essential spectrum and 
absolutely continuous spectrum, respectively, and 
${\mathcal H}_{ac}(A)$ and ${\mathcal H}_{pp}(A)$ are the absolutely continuous subspace 
for $A$
and  the closure of the linear hull of eigenvectors for $A$, respectively. Generic 
points on ${\mathcal M}$
are denoted by $p, \dots,$ or $X, Y, \dots,$ while those in the ends ${\mathcal M}_j$ 
are often written
as $(x, y)$. ${\bf N}$ denotes the set of all positive integers. When $\bf h$, $I \subset \R$ 
is an interval and $d \mu$ is a measure on $I$, $L^2(I, {\bf h}; d\mu)$ denotes the space of
all $\bf h$-valued $L^2-$functions on $I$ with respect to $d\mu$.


\section{2-dimensional hyperbolic orbifolds and conical singularties}


\subsection{Fuchsian groups}
The upper-half space model of 2-dimensional hyperbolic space ${\bf H}^2$ is  ${\bf C}_+ = \{z = x + iy \, ; y > 0\}$ equipped with the metric
\begin{equation}
ds^2 = \frac{(dx)^2 + (dy)^2}{y^2}.
\label{S2metric}
\end{equation}
The infinity of ${\bf H}^2$ is 
$$
\partial{\bf C_+} = 
{\bf R}\cup{\infty}.
$$
${\bf H}^2$ admits an action of $SL(2,{\bf R})$ defined by
\begin{equation}
SL(2,{\bf R})\times{\bf C}_+ \ni (\gamma,z) \to \gamma\cdot z:= \frac{az + b}{cz + d}, \quad \gamma = \left(
\begin{array}{cc}
a & b \\
c & d
\end{array}
\right).
\label{S2Action}
\end{equation}
The right-hand side, M{\"o}bius transformation, is an isometry on ${\bf H}^2$. The 
mapping : $\gamma \to \gamma\cdot$ is 2 to 1, and the corresponding factor group of 
M{\"o}bius transformations is isomorphic to 
$PSL(2,{\bf R}) = SL(2,{\bf R})/\{\pm I\}$.
For $\gamma \neq \pm I$, the transformation (\ref{S2Action}) is classified into 3 categories :
\begin{equation}
\begin{split}
elliptic &\Longleftrightarrow {\rm there \ is \ only \ one \ fixed \ point \ in }\ {\bf C}_+ 
\\
&\Longleftrightarrow |{\rm tr}\, \gamma| < 2, \\
parabolic &\Longleftrightarrow {\rm there \ is \ only \ one \ degenerate \ fixed \ point \ on} \ \partial{\bf C}_+ \\
&\Longleftrightarrow |{\rm tr}\, \gamma| = 2, \\
hyperbolic &\Longleftrightarrow {\rm there \ are \ two \  fixed \ points \ on} \ \partial{\bf C}_+ \\
&\Longleftrightarrow |{\rm tr}\, \gamma| > 2.
\end{split}
\nonumber
\end{equation}

Let $\Gamma$ be a discrete subgroup, {\it Fuchsian group}, of $SL(2,{\bf R})$, and 
$\mathcal M = \Gamma\backslash{\bf H}^2$  by the action 
(\ref{S2Action}). $\Gamma$ is said to be {\it geometrically finite} if the fundamental domain
$\Gamma\backslash{\bf H}^2$ is 
chosen to be a finite-sided convex polygon. The sides are then geodesics of ${\bf H}^2$. 
The geometric finiteness is equivalent to that $\Gamma$ is finitely generated (\cite{Ka92}, 
p. 104). Let us give two simple but important examples. 

\subsubsection{Parabolic cyclic group}
Consider the cyclic group $\Gamma$ generated by the action $z \to z + \tau$. This is parabolic 
with fixed point $\infty$. The associated fundamental domain is then
$[-\tau/2,\tau/2]\times(0,\infty)$ with the sides $x = \pm \tau/2$ being geodesics.
The Riemann surface $\mathcal M$ is then equal to $S^{\tau/2 \pi} \times (0, \infty)$,
which is a hyperbolic manifold with metric (\ref{S2metric}). 
It has two infinities : $S^{\tau/2 \pi}\times \{0\}$ and $\infty$. The part 
$S^{\tau/2 \pi}\times(0,1)$ has an infinite volume. 
The part $S^{\tau/2 \pi}\times(1,\infty)$ has a finite volume, and is called the {\it cusp}.

\subsubsection{Hyperbolic cyclic group}
Another simple example is the cyclic group generated by the hyperbolic action $z \to \lambda z$, $\lambda > 1$. The sides of the fundamental domain 
$ \{1 \leq |z| \leq \lambda\}$ are semi-circles orthogonal to $\{y = 0\}$, which are geodesics. 
The quotient manifold is diffeomorphic to $S^{(\log\lambda)/2\pi}\times(-\infty,\infty)$. It is parametrized by $(t,r)$, where $t \in {\bf R}/(\log\lambda){\bf Z}$ and $r$ is the signed distance from the segment 
$\{(0,s)\, ; 1 \leq s \leq \lambda\}$. The metric is then written as
\begin{equation}
ds^2 = (dr)^2 + \cosh^2r\,(dt)^2.
\label{S2Funnelmetric}
\end{equation}
The part $r > 0$ (or $r < 0)$ is called the {\it funnel}. Letting $y = 2e^{-r},\, r>0$, one can rewrite (\ref{S2Funnelmetric}) as
\begin{equation}
ds^2 = \Big(\frac{dy}{y}\Big)^2 + \Big(\frac{1}{y} + \frac{y}{4}\Big)^2(dt)^2.
\nonumber
\end{equation}
Therefore, the funnel is regarded as a perturbation of the infinite volume part 
$S^{(\log\lambda)/2 \pi}\times(0,1)$ of the fundamental domain for the parabolic cyclic group.


\subsection{Classification of 2-dimensional hyperbolic manifolds}
The set of limit points of a Fuchsian group $\Gamma$, denoted by $\Lambda(\Gamma)$, is defined as follows : $w \in \Lambda(\Gamma)$ if there exist $z_0 \in {\bf C}_+$ and 
$\gamma_n \in \Gamma, \, \gamma_n \neq \hbox{I},$ such that $\gamma_n\cdot z_0 \to w$.
 Since $\Gamma$ acts discontinuously on ${\bf C}_+$, $\Lambda(\Gamma) \subset \partial{\bf H}^2 = \partial{\bf C}_+$. There are only 3 possibilities.
\begin{itemize}
\item 
({\it Elementary}) : 
$\Lambda(\Gamma)$ is a finite set.

\item
({\it The 1st kind}) :
$\ \Lambda(\Gamma) = \partial{\bf H}^2$.

\item
({\it The 2nd kind}) :
$\ \Lambda(\Gamma)$ is a perfect (i.e. every point in $\Lambda(\Gamma)$ is an accumulation point of $\Lambda(\Gamma)$), nowhere dense set of $\partial{\bf H}^2$.
\end{itemize}
 
Any elementary group is either cyclic or is conjugate in $PSL(2,{\bf R})$ to a group generated by $\gamma\cdot z = \lambda z$, $(\lambda > 1)$, and $\gamma'\cdot z = - 1/z$ (see \cite{Ka92}, Theorem 2.4.3).

For non-elementary case, we have the following theorem (\cite{Bor07}, Theorem 2.13). 
Although \cite{Bor07} deals with the case without elliptic fixed points, this theorem holds for the case with elliptic fixed points.


\begin{theorem}
Let ${\mathcal M} = \Gamma\backslash{\bf H}^2$ be a non-elementary geometrically finite hyperbolic manifold. Then there exists a compact subset $\mathcal K$ such that $\mathcal M \setminus {\mathcal K}$ is a finite disjoint union of cusps and funnels.
\end{theorem}

The regions mentioned above, i.e. fundamental domains of parabolic cyclic groups, hyperbolic cyclic groups, and non-elementary geometrically finite groups are the models of hyperbolic spaces to be dealt with in this paper.

Other important theorems are the following (see \cite{Ka92}, Theorems 4.5.1, 4.5.2 and 4.1.1).


\begin{theorem}
A Fuchsian group is of the 1st kind if and only if its fundamental domain has a finite area.
\end{theorem}


\begin{theorem}
A Fuchsian group of the 1st kind is geometrically finite.
\end{theorem}

 For the Fuchsian group of the 1st kind, therefore, the ends of its fundamental domain are always cusps. In this case, usually it is compactified around parabolic fixed points and made to a compact Riemann surface. The automorphic functions associated with this group turn out to be algebraic functions on this Riemann surface (see \cite{Leh66}). 
 
 It is well-known that there is a 1 to 1 correspondence between the compact Riemann surfaces and the fields of algebraic functions. 
This suggests a general idea that a surface will be determined by a set of functions on it. What we pursue in this paper is an analogue of this fact. Asymptotically hyperbolic manifolds, more generally non-compact Riemannian manifolds with good structure at infinity will be determined by the set of solutions to the Helmholtz equation, more precisely, by the asymptotic behavior at infinity of solutions to the Helmholtz equation.
 Before going into the detail of this issue, we need to recall the notion of orbifolds.



\subsection{Elliptic fixed point and analytic structure of $\Gamma\backslash{\bf H}^2$.}

Now, we study the analytic structure of $\mathcal M = \Gamma\backslash{\bf H}^2$, where $\Gamma$ is a Fuchsian group.
 Let $\mathcal M_{sing}$ be the set of all elliptic fixed points in $\mathcal M $. Under the assumption of geometric finiteness, 
$\mathcal M_{sing}$ is a finite set. 


\begin{lemma}
Let $p \in \mathcal M_{sing}$, and 
$$
\mathcal I(p) = \{\gamma \in \Gamma\, ; \, \gamma\cdot p = p \}
$$  
the isotropy group of $p$. Then, it is a finite cyclic group, and its generator $\gamma_0$ satisfies
\begin{equation}
\frac{w - p}{w - \overline{p}} = e^{2\pi i/n}\frac{z - p}{z - \overline{p}}, \quad w = \gamma_0\cdot z
\label{S2ztow}
\end{equation}
for some $n=n_p \in {\bf N}$.
\end{lemma}

\noindent
{\bf Proof.} Recall that the cross ratio
$$
(z_1,z_2,z_3,z_4) = \frac{z_1 - z_3}{z_1 - z_4}\cdot\frac{z_2 - z_4}{z_2 - z_3}
$$
is invariant by the fractional linear transformation 
$z \to w = \gamma\cdot z= (az + b)(cz + d)^{-1}$. 
Suppose $p, q \in {\bf C}$ are the fixed points of $z \to w$. Then, since 
$z = \infty$ 
is mapped to $w = a/c$, we have
$(w,a/c,p,q) = (z,\infty,p,q)$,
which implies
$$
\frac{w-p}{w-q} = \kappa\frac{z-p}{z-q}, \quad \kappa = \frac{a-cp}{a-cq}.
$$
For the elliptic case, $|\kappa| = 1$, since $q = \overline{p}$.
By the linear fractional transformation 
$T(z) = (z - p)/(z - \overline{p})$, $\gamma$
  is written as $\gamma = T^{-1}\kappa T$.
Therefore, $\mathcal I(p)$ is isomorphic to a discrete subgroup of $SO(2)$, which proves (\ref{S2ztow}). \qed

\medskip
{\nntext To introduce the analytic structure near $p$, we let}
 $\iota$ be the canonical projection 
\begin{equation}
\iota : {\bf H}^2 \ni z \to [z] = \{g\cdot z\, ; \, g \in \Gamma\} \in \Gamma\backslash{\bf H}^2.
\nonumber
\end{equation}
Using $n$ from (\ref{S2ztow}),
we introduce the local coordinates $\varphi_p(\iota(z))$ near $p$ by
\begin{equation}
\zeta := \varphi_p(\iota(z)) = T(z)^n= \left( \frac{z-p}{z- {\overline p}} \right)^n,
\quad \zeta(p)=0.
\nonumber
\end{equation}
 Identifying $z$ and $\iota(z)$, we have as $\zeta \to 0$
$$
z = \frac{p - \overline{p}\zeta^{1/n}}{1 - \zeta^{1/n}} = p + (p-{\overline p})\zeta^{1/n} + \cdots.
$$
Therefore,
\begin{equation}
\frac{(dx)^2 + (dy)^2}{y^2} = \frac{dz\,d\overline{z}}{({\rm Im}\,z)^2} = 
\frac{\big|dz/d\zeta\big|^2}{({\rm Im}\,z)^2}d\zeta d\overline{\zeta}.
\label{S2dzetadzetabar}
\end{equation}
Direct computation entails
$$
\frac{dz}{d\zeta} = \frac{p -\overline{p}}{n}\zeta^{1/n-1}(1 - \zeta^{1/n})^{-2},
$$
$$
{\rm Im}\,z = \frac{p-\overline{p}}{2i}\frac{1 - |\zeta^{1/n}|^2}{|1-\zeta^{1/n}|^2}.
$$
Therefore, we have
\begin{equation}
\frac{|dz/d\zeta|^2}{({\rm Im} \, z)^2} = \frac{4}{n^2}|\zeta|^{-\lambda}\left(1 - |\zeta|^{2/n}\right)^{-2}, \quad \lambda = 2 - \frac{2}{n}.
\label{S2dzdzeta}
\end{equation}
Note that $1 \leq  \lambda < 2$.
The volume element and the Laplace-Beltrami operator are then rewritten as
\begin{equation}
\frac{dx\wedge dy}{y^2} = \frac{i}{2y^2}dz\wedge d\overline{z} = \frac{i\big|dz/d\zeta\big|^2}{2({\rm Im}\,z)^2}d\zeta\wedge d\overline{\zeta},
\nonumber
\end{equation}
\begin{equation}
y^2\big(\partial_x^2 + \partial_y^2\big) = 4({\rm Im}\,z)^2\frac{\partial^2}{\partial z\partial\overline{z}} = 
\frac{4({\rm Im}\,z)^2}{\big|dz/d\zeta|^2}\frac{\partial^2}{\partial\zeta\partial\overline{\zeta}}.
\nonumber
\end{equation}
Both of them have singularities at $p$. However, if $f, g$ are $C^{\infty}$-functions 
with respect to $\zeta$ supported 
near $p$, we have
\begin{equation}
\int_{\mathcal M}y^2\big(\partial_x^2 + \partial_y^2\big)f\cdot {\overline g}\,\frac{dxdy}{y^2} 
= 2i
\int_{|\zeta|<\epsilon} 
\frac{\partial^2}{\partial\zeta\partial\overline{\zeta}}f\cdot {\overline g}\, 
d\zeta d\overline{\zeta}
\\ \nonumber
=-2i \int_{|\zeta|<\epsilon} \frac{\partial f}{\partial\zeta}\,
{\overline{\frac{\partial g}{\partial\zeta}}}\,d\zeta d\overline{\zeta}.
\nonumber
\end{equation}
What is important is that the singularity of the volume element and that of the Laplace-Beltrami operator cancel. 

Let $\mathcal M_{sing} = \{p_1, \cdots, p_L\}$. We take a small open set 
$U_j \subset \mathcal M$ such that $p_j \in U_j$, $U_i\cap U_j = \emptyset$ if $i \neq j$. 
We construct a smooth partition of unity $\{\chi_j\}_{j=0}^L$ such that 
${\rm supp}\,\chi_j \subset U_j$, $j = 1,\cdots,L$, and $\sum_{j=0}^L\chi_j = 1$ on $\mathcal M$. 
We put
\begin{equation}
dV_H^{(j)} = \left\{
\begin{split}
& \frac{dx\wedge dy}{y^2}=\frac{i}{2} \frac{ dz\wedge d\overline{z}}{({\rm Im}\, z)^2} 
\quad (j = 0), \\
& \frac{i\big|dz/d\zeta\big|^2}{2({\rm Im}\, z)^2}d\zeta\wedge d\overline{\zeta} \quad (j \neq 0),
\end{split}
\right.
\label{S2dVHj}
\end{equation}
\begin{equation}
dV_E^{(j)} = \left\{
\begin{split}
& \frac{i}{2} dz \wedge d\overline{z} \quad (j = 0), \\
& \frac{i}{2}d\zeta\wedge d\overline{\zeta} \quad (j \neq 0),
\end{split}
\right.
\label{S2dVEj}
\end{equation}
and define a quadratic form $Q_{AS}[u,v]$ by
\begin{equation}
Q_{AS}[u,v] = \sum_{j=0}^L\int_{\mathcal M}\chi_j u\,\overline{v}\ dV_H^{(j)} + 
\sum_{j=0}^L\int_{\mathcal M}\chi_i\nabla u\cdot\nabla\overline{v}\ dV_E^{(j)},
\nonumber
\end{equation}
where
\begin{equation}
\nabla = \left\{
\begin{split}
& (\partial_x,\partial_y) \quad (j = 0), \\
& (\partial_t,\partial_s) \quad (j \neq 0), \quad (\zeta = t + is).
\end{split}
\right.
\nonumber
\end{equation}

Let $L^2(\mathcal M)$ be the Hilbert space of $L^2$-functions on $\mathcal M$ with respect to the measure $dxdy/y^2$. 
As is easily seen, $\sqrt{Q_{AS}(u,u)}$ defines a norm on $C_0^{\infty}(\mathcal M\setminus\mathcal M_{sing})$.
Let $D(Q_{AS})$ be the completion of $C_0^{\infty}(\mathcal M\setminus\mathcal M_{sing})$ with respect to the norm $\sqrt{Q_{AS}[u,u]}$. This is the counterpart of the 1st-order Sobolev space on $\mathcal M$.


\begin{lemma} Let $\Gamma$ be a geometrically finite Fuchsian group.  Then,
for any compact set $K \subset \Gamma\backslash{\bf H}^2$, the imbedding
$$
D(Q_{AS}) \ni u \to u\big|_K \in L^2(K)
$$
is compact.
\end{lemma}

\noindent
{\bf Proof.} This is obvious if $K$ does not contain elliptic fixed points. Around an elliptic fix point $p_j$ $(1 \leq j \leq L)$, we take local coordinate $\zeta = t + is$ as above, and for a suffiently small $r > 0$, let $B_r = \{(t,s)\,;\,t^2 + s^2 < r^2\}$. Then, by (\ref{S2dVEj}), if $u \in D(Q_{AS})$ has a support in $B_r$,
\begin{equation}
\int_{B_r}|u|^2dtds \leq C\int_{B_r}|u|^2dV_E^{(j)},
\label{S2BrdtdsBrdVH}
\end{equation}
with a constant $C > 0$. By the Sobolev imbedding $H^s({\bf R}^d) \subset L^p_{loc}({\bf R}^d)$, where $0 \leq s < d/2$, $p = 2d/(d-2s)$, we have
\begin{equation}
H^1({\bf R}^2) \subset L^p_{loc}({\bf R}^2), \quad \forall p > 2,
\label{S2Sobolev}
\end{equation}
with continuous inclusion. 

We take $\alpha, \beta$ such that $\alpha^{-1} + \beta^{-1} = 1$, $1 < \alpha < 2/\lambda$. Then  by (\ref{S2dzdzeta}), (\ref{S2dVHj}), and  H{\"o}lder's inequality
$$
\int_{B_{\delta}}|u|^2dV_H^{(j)} \leq C\int_{B_{\delta}}r^{-\lambda}|u|^2dtds
\leq C\left(\int_{B_{\delta}}r^{-\lambda\alpha}dtds\right)^{1/\alpha}
\left(\int_{B_{\delta}}|u|^{2\beta}dtds\right)^{1/\beta},
$$
where $r = (s^2 + t^2)^{1/2}$.
Since $\lambda\alpha < 2$, the 1st term of the most right-hand side tends to 0 when $\delta \to 0$. To the 2nd term of the most right-hand side we apply (\ref{S2Sobolev}). Using (\ref{S2BrdtdsBrdVH}), for any $\epsilon > 0$ there exists $\delta > 0$ such that
$$
\int_{B_{\delta}}|u|^2dV_H^{(j)} \leq \epsilon\left(\int_{B_{\delta}}|u|^2dV_H^{(j)} + \int_{B_{2\delta}}|\nabla u|^2dV_E^{(j)}\right).
$$
Suppose we are given a bouded sequence $\{u_n\}$ in $D(Q_{AS})$. Then the integral of $|u_n|^2$ over $B_{\delta}$ with respect to the measure $dV_H^{(j)}$ can be made small uniformly in $n$. Outside $B_{\delta}$, we use the usual Rellich theorem. This proves the lemma. \qed

\medskip
Let $H_{AS}$ be  the Laplce-Beltrami operator $- \Delta_g-1/4$ on 
${\mathcal M }$,  defined through the quadratic form $Q_{AS}[u,v]$. 
It is well-known that
\begin{equation}
D(H_{AS}) \subset D(H_{AS}^{1/2}) = D(Q_{AS}).
\label{S2DHorb}
\end{equation} 


\begin{cor}
 $\chi (H_{AS} - z)^{-1},\, z \notin {\bf R}$, is compact on $L^2(\Gamma\backslash{\bf H}^2)$ for any $\chi \in C_0^{\infty}(\Gamma\backslash{\bf H}^2)$.
\end{cor}

\noindent
{\bf Proof.} This follows from Lemma 2.5 and (\ref{S2DHorb}). \qed

\medskip
Using these facts, one can discuss the forward problem, i.e. the spectral theory for $H_{AS}$ 
in the same way as in \cite{IsKu10}. In order to discuss the inverse problem, however, it is 
more appropriate to change the differentiable structure around singular points and introduce 
the notion of conical singularities.


\subsection{Manifolds with conical singularities. Orbifolds.} 
Let us repeat the definition of a Riemannian surface with conical singularities.
We warn the reader not to confuse it with the Riemann surface (the 1-dimensional complex manifold).


\begin{definition}
{\it A $C^{\infty}$-surface $\mathcal M$ is said to be a 
Riemannian surface with conical singularities if there exists a discrete 
subset $\mathcal M_{sing}$ of $\mathcal M$ such that \\
\noindent 
(i) there exists a smooth Riemannian metric $g$ on $\mathcal M\setminus\mathcal M_{sing}$, \\
\noindent
(ii) for each $p \in \mathcal M_{sing}$, there is an open neighborhood $U_p$ of $p$ such that the assumption (A-2) is satisfied on $U_p$.}
\end{definition}

Let $\mathcal M$ be a Riemannian surface with conical singularities. Then, near $p \in \mathcal M_{sing}$, letting $x^1 = r\cos\theta$, $x^2 = r\sin\theta$, we see that the metric $g = g_{ij}dx^idx^j$ satisfies
\begin{equation}
C^{-1}I \leq \left(g_{ij}\right) \leq CI, \quad C>1.
\label{S2boundedmetric}
\end{equation}
This shows that, although the metric $g$ may be singular at $\mathcal M_{sing}$, the $H^1$-norm
\begin{equation}
\|u\|_{H^1(\mathcal M)} = \left(\int_M|u|^2\sqrt{g}dx + 
\int_Mg^{ij}\frac{\partial u}{\partial x^i}\frac{\partial\overline u}{\partial x^j}\sqrt{g}dx\right)^{1/2}
\label{S2H1normconial}
\end{equation}
can be introduced in the same way as in the case of $C^{\infty}$-Riemannian manifold. In particular, the following lemma holds. 


\begin{lemma}
Let $\mathcal M$ be a surface satisfying (A-1) $\sim$ (A-4), and $\Delta_g$ its Laplace-Beltrami 
operator. Then $- \Delta_g - 1/4$ has a self-adjoint realization through the quadratic form, which is denoted by $H$. Then, for any 
$\chi \in C_0^{\infty}(\mathcal M)$, $\chi(H - z)^{-1}$, $z \not\in {\bf R}$, is a compact 
operator on $L^2(\mathcal M)$.
\end{lemma}

Next we return to  $\mathcal M = \Gamma\backslash{\bf H}^2$,
where $\Gamma$ is a Fuchsian group.
We show that (A-2) is satisfied around an elliptic fixed point $p \in \mathcal M$. By (\ref{S2dzdzeta}), putting $\zeta = \rho e^{i\theta}$, the metric (\ref{S2dzetadzetabar}) takes the form
$$
\frac{4}{n^2}|\zeta|^{-\lambda}(1-|\zeta|^{2/n})^{-2}d\zeta d\overline{\zeta} = (1 - \rho^{2/n})^{-2}\frac{4}{n^2}\rho^{-\lambda}\left((d\rho)^2 + \rho^2(d\theta)^2\right).
$$
Putting $t = 2\rho^{1-\lambda/2} = 2\rho^{1/n}$, we can rewrite it as
$$
(1-t^2/4)^{-2}\left((dt)^2 + \frac{t^2}{n^2}(d\theta)^2\right).
$$
Solving $dr = (1 - t^2/4)^{-1}dt$, we have
$$
r = \log\frac{2+t}{2-t} = \log\frac{1+\rho^{1/n}}{1-\rho^{1/n}}.
$$
Therefore $\rho^{1/n} = (e^r-1)/(e^r+1)$, and (\ref{S2dzetadzetabar}) takes the form
\begin{equation}
\frac{(dx)^2 + (dy)^2}{y^2} = (dr)^2 + \frac{1}{n^2}\sinh^2r(d\theta)^2.
\label{S2metricconical}
\end{equation}
This shows that (A-2) is satisfied for any $p \in \mathcal M_{sing}$. We cover 
$\mathcal M\setminus\mathcal M_{sing}$ by  standard local coordinate patches of the 
quotient 
Riemannian surface $\Gamma\backslash{\bf H}^2$. Therefore,
 $\mathcal M$ is a Riemannian surface with conical singularties.
{\nntext Actually, the structure of conical singularities on 
$\mathcal M =\Gamma\backslash{\bf H}^2$
is of a special form, making it a {it Riemannian orbifold}.

To define an (orientable) 2D-Riemannian orbifold, let 
  $\mathcal M$ be a 2D-manifold. 
  Suppose there exists a discrete subset $\mathcal M_{sing} \subset \mathcal M$ 
  such that  
  $\mathcal M \setminus\mathcal M_{sing}$ is an orientable Riemannian manifold 
  with a $C^{\infty}$-Riemannian metric $g$.}
 We assume that 
 each point  $p \in \mathcal M_{sing}$
has a  neighborhoods for which the following properties  hold (see \cite{Sa57}, \cite{Thurs});
 \begin{enumerate}
 \item[{\bf (B-1)}]  There exists an open set $\widetilde U_p^\e$ in ${\bf R}^2$, containing the 
 origin $0$
 and equipped with a Riemannian metric $\widetilde g_p$,
such that,  with respect to $\widetilde g_p$, $\widetilde U_p^\e$ is the ball of 
radius $\e$ 
centered at $0$.
\item[{\bf (B-2)}]  There is a finite group of rotations $\Gamma_{n_p} \subset SO(2)$ of order 
$ n_p >1,$ 
so that $\widetilde g_p$
is invariant with respect to the action of $\Gamma_{n_p}$.
\item[{\bf (B-3)}]  $U_p^\e \sim \widetilde U_p^\e /\Gamma_{n_p}$, where $U_p^\e$ is  the ball of radius $\e$ 
on ${\mathcal M}$,
centered at $p$, and $\sim$ stands for the isometry.
 \end{enumerate}


If these assumptions are satisfied, we say that $\mathcal M$ is a 2-dimensional  Riemannian orbifold.
We call $n_p$ the {\it order} of $p \in \mathcal M_{sing}$.
 For the neighborhoods defined in condition B-3 we denote by $\pi_p:\widetilde U_p^\e
\to  U_p^\e$ the associated canonical projections and say that $(\widetilde U_p^\e,\widetilde g_p)$
is the uniformizing  cover of $(U_p^\e, g_p)$.

 A homeomorphism $\Phi$ between Riemannian orbifolds $\mathcal M^{(1)}$ and $\mathcal M^{(2)}$ is said to be an {\it orbifold isomorphism} if it has the following properties: \\
\noindent
(1) $\Phi : \mathcal M^{(1)}\setminus\mathcal M^{(1)}_{sing} \to \mathcal M^{(2)}\setminus\mathcal M^{(2)}_{sing}$ is a Riemannian isometry. \\
\noindent
(2) For any $p^{(1)} \in \mathcal M^{(1)}_{sing}$ and $p^{(2)} = \Phi(p^{(1)})$, $\Phi : U_{p^{(1)}}^{\epsilon} \to U_{p^{(2)}}^{\epsilon}$ is lifted to an isometry between the coverings $\widetilde \Phi :  \widetilde U_{p^{(1)}}^{\epsilon} \to \widetilde U_{p^{(2)}}^{\epsilon}$.

To bridge the notion of a surface with conical singularities with that of a 2-dimensional Riemannian orbifold, 
note that an orbifold singularity is a particular case of a conical singularity characterized by two properties: 
\begin{condition} \label{cond_orbifold} 
 {\bf i.} $C_p= \big(1/n_p\big)^2$.

 {\bf ii.} The metric
tensor (\ref{S1metricaroundsingular}), rewritten in  coordinates 
$x^1=r \cos(\theta/n_p),$ $ x^2=r \sin(\theta/n_p)$, being continued periodically onto 
${\tilde U}_\e(0)=\{r < \e  \}$ is smooth.  
\end{condition}

{\nntext Returning to $\Gamma\backslash{\bf H}^2$ and using equation 
(\ref{S2metricconical}), straighforward 
calculations show  that each singular point $p \in \Gamma\backslash{\bf H}^2$ satisfies
conditions {\bf i., ii.}
}

\medskip

Let us summarize what we have done in this section. For the Fuchsian group $\Gamma \in SL(2,{\bf R})$,  $\mathcal M = \Gamma\backslash{\bf H}^2$
{\nntext has a structure of a 2D-Riemannian orbifold.} It is a Riemann surface, i.e. 1-dimensional complex manifold without singularities.  By changing the differentiable structure around $\mathcal M_{sing}$ = the set of the elliptic fixed points, $\mathcal M$ is regarded as a Riemannian surface with conical singularities. These two local coordinate systems have the following features.
\begin{itemize}
\item They coincide except for a small neighborhood of $\mathcal M_{sing}$, and give an equivalent $C^{\infty}$-differentiable structure on $\mathcal M \setminus{\mathcal M}_{sing}$.
\item They equip $\mathcal M \setminus{\mathcal M}_{sing}$ with the hyperbolic metric, which is singular at $\mathcal M_{sing}$ in the case of orbifold.
\item The associated Laplace-Beltarmi operators are unitarily equivalent.
\end{itemize}

It follows from these properties that the associated (generalized) S-matrices coincide, since they are defined by the asymptotic behavior at infinity of solutions to the Helmholtz equations.

This new coordinate system resolves the singularities of the hyperbolic metric at elliptic fixed points, which makes the proof of local compactness of the resolvent easier. The merit of introducing the notion of conical singularities is not restricted here, however. It is used effectively in the inverse problem in \S 5. On the other hand, the original coordinate system is analytic even at elliptic fixed points. This fact will be used in \S 6 to discuss the orbifold isomorphism.


\section{Spectral theory for asymptotically hypebolic manifolds}
In \cite{IsKu10}, for manifolds without conical singularities,
we have already studied spectral properties of the Laplace-Beltrami operators on asymptotically 
hyperbolic manifolds : limiting absorption principle for the resolvent, spectral 
representations, S-matrices. Thanks to Lemma 2.8, and also to the fact that $\mathcal M_{sing}$ is 
a finite set, the proof of the above facts works well without any change. We shall explain below 
the basic ideas for 
this forward problem and summarize the results.

Let $A$ be a self-adjoint operator in a Hilbert space $\mathcal H$. If $\lambda \in \sigma(A)$, the limit $\lim_{\epsilon\to0}(A - \lambda \mp i\epsilon)^{-1}$ does not exist in ${\bf B}(\mathcal H;\mathcal H)$. However, in some important cases, when $\lambda \in \sigma_{c}(A)$, 
it is possible to define $\lim_{\epsilon\to0}(A - \lambda \mp i\epsilon)^{-1}$. This is achieved
by choosing suitable Banach spaces $\mathcal H_+, \mathcal H_-$ 
satisfying
$$
\mathcal H_+ \subset \mathcal H \subset \mathcal H_-, 
$$
with continuous injections, so that
\ba
\lim_{\e \to 0} \left(A-\la \mp i\e  \right)^{-1} \in {\bf B}(\mathcal H_+;\mathcal H_-).
\ea  
This fact is usually called the limiting absorption principle.
For $A = - \Delta$ in ${\bf R}^n$, the best choice of $\mathcal H_{\pm}$ are the 
Besov type spaces $\mathcal B, \mathcal B^{\ast}$ found by Agmon-H{\"o}rmander 
\cite{AgHo76}.
We first define a counterpart of $\mathcal B, \mathcal B^{\ast}$ in the case of hyperbolic spaces.


\subsection{Besov type spaces}
Let $\bf h$ be a Hilbert space endowed with inner product $(\;,\;)_{\bf h}$ and norm $\|\cdot\|_{\bf h}$.
We decompose $(0,\infty)$ into $(0,\infty) = \cup_{k\in {\bf Z}}I_k$, where
\begin{equation}
I_k = \left\{
\begin{array}{cc}
\big(\exp(e^{k-1}),\exp(e^k)\big], & k \geq 1, \\
\big(e^{-1},e\big], & k = 0, \\
\big(\exp(-e^{|k|}),\exp(-e^{|k|-1})\big], & k \leq -1. 
\end{array}
\right.
\nonumber
\end{equation}
Let $\mathcal B = \mathcal B({\bf h})$ be the Banach space of ${\bf h}$-valued function on $(0,\infty)$ such that
\begin{equation}
\|f\|_{\mathcal B} = \sum_{k\in{\bf Z}}e^{|k|/2}
\left(\int_{I_k}\|f(y)\|_{\bf h}^2\, \frac{dy}{y^2}\right)^{1/2} < \infty.
\label{S3Bspace}
\end{equation}
The dual space of $\mathcal B$ is identified with the space equipped with norm 
\begin{equation}
\|u\|_{\mathcal B^{\ast}} = \left(\sup_{R>e}\frac{1}{\log R}\int_{1/R}^R\|u(y)\|_{\bf h}^2\, \frac{dy}{y^2}\right)^{1/2} < \infty.
\label{S3Bast}
\end{equation}
For example, for $\phi \in {\bf h}$, $y^{1/2}\phi$ belongs to $\mathcal B^{\ast}$.
We also use the following weighted $L^2$-space: for $s \in {\bf R}$, 
\begin{equation}
L^{2,s} \ni u \Longleftrightarrow \|u\|_{s} = \left(\int_0^{\infty}
(1 + |\log y|)^{2s}\|u(y)\|_{\bf h}^2\,\frac{dy}{y^2}\right)^{1/2} < \infty.
\label{S3L2s}
\end{equation}
For $s > 1/2$, the following inclusions hold:
\begin{equation}
L^{2,s} \subset \mathcal B \subset L^{2,1/2} \subset L^2 \subset L^{2,-1/2} \subset \mathcal B^{\ast} \subset L^{2,-s}.
\label{S3Inclusion}
\end{equation}

If $u, v \in \mathcal B^{\ast}$ satisfy
\begin{equation} \label{S3Bast0}
 \lim_{R\to\infty}\frac{1}{\log R}\int_{1/R}^{1/2}\|u(y) - v(y)\|^2_{\bf h}\frac{dy}{y^2} = 0,
\quad \lim_{R\to\infty}\frac{1}{\log R}\int_{2}^R\|u(y) - v(y)\|^2_{\bf h}\frac{dy}{y^2} = 0,
\end{equation}
we regard that $u$ and $v$ have the same asymptotic behavior at infinities,  $y = 0$, and 
$y = \infty$, correspondingly. We have the following lemma.


\begin{lemma} 
For $u \in \mathcal B^{\ast}$, the following two assertions are equivalent.
\begin{equation}
\lim_{R\to\infty}\frac{1}{\log R}\int_{1/R}^R\|u(y)\|_{\bf h}^2\,\frac{dy}{y^2} = 0.
\end{equation}
\begin{equation}
\lim_{R\to\infty}\frac{1}{\log R}\int_0^{\infty}\rho\Big(\frac{\log y}{\log R}\Big)\|u(y)\|_{\bf h}^2\,\frac{dy}{y^2} = 0, \quad \forall \rho \in C_0^{\infty}((0,\infty)).
\end{equation}
\end{lemma}

The proof of the above results are given in \cite{IsKu10}, Chap. 1, \S 2. 


\subsection{Bessel functions} 
We use the following knowledge of Bessel functions.
 For the details, see \cite{Wa62}.
The modified Bessel function (of the 1st kind) $I_{\nu}(z)$, with parameter 
$\nu \in {\bf C}$, 
is defined by
\begin{equation}
I_{\nu}(z) = \left(\frac{z}{2}\right)^{\nu}\sum_{n=0}^{\infty}
\frac{(z^2/4)^n}{n!\,\Gamma(\nu + n + 1)}, \quad
z \in {\bf C}\setminus(-\infty,0].
\label{S3DefinitionofInu}
\end{equation} 
It is related to the Bessel function $J_{\nu}(z)$ as follows 
\begin{equation}
I_{\nu}(y) = e^{-\nu\pi i/2}J_{\nu}(iy), \quad y > 0.
\nonumber
\end{equation}
The following function $K_{\nu}(z)$ is also called the modified Bessel function, or the K-Bessel function, or sometimes the Macdonald function:
\begin{equation}
K_{\nu}(z) = \frac{\pi}{2}\frac{I_{-\nu}(z) - I_{\nu}(z)}{\sin(\nu\pi)}, \quad
\nu \notin {\bf Z},
\label{S3KnuandInu}
\end{equation}
\begin{equation}
K_n(z) = K_{-n}(z) = \lim_{\nu \to n}K_{\nu}(z), \quad 
n \in {\bf Z}.
\nonumber
\end{equation}
These $I_{\nu}(z), K_{\nu}(z)$ solve the following equation 
\begin{equation}
z^2u'' + zu' - (z^2 + \nu^2)u = 0,
\label{S3BesselDiffEq}
\end{equation}
and have the following asymptotic expansions as $|z| \to \infty$:
\begin{equation}
I_{\nu}(z) \sim \frac{e^z}{\sqrt{2\pi z}} +
 \frac{e^{-z + (\nu + 1/2)\pi i}}{\sqrt{2\pi z}}, \quad |z| \to \infty, \quad
- \frac{\pi}{2} < {\rm arg}\,z < \frac{\pi}{2},
\label{S3Inunearinfty}
\end{equation}
\begin{equation}
K_{\nu}(z) \sim \sqrt{\frac{\pi}{2z}}e^{-z}, \quad 
|z| \to \infty, \quad - \pi < {\rm arg}\,z < \pi.
\label{S3Knunearinfty}
\end{equation}
The asymptotics as $z \to 0$ are as follows:
\begin{equation}
I_{\nu}(z) \sim \frac{1}{\Gamma(\nu + 1)}\left(\frac{z}{2}\right)^{\nu},
\label{S3Inunear0}
\end{equation}
\begin{equation}
K_{\nu}(z) \sim \frac{\pi}{2\sin(\nu\pi)}
\left(\frac{1}{\Gamma(1 - \nu)}\left(\frac{z}{2}\right)^{-\nu} - 
\frac{1}{\Gamma(1 + \nu)}\left(\frac{z}{2}\right)^{\nu}\right), \quad 
\nu \not\in {\bf Z}
\label{S3Knunear0}
\end{equation}
\begin{equation}
K_n(z) \sim \left\{
\begin{split}
& - \log z, \quad n = 0, \\
& 2^{n-1}(n - 1)!z^{-n}, \quad n = 0, 1, 2, \dots.
\end{split}
\right.
\nonumber
\end{equation}


\subsection{Spectral properties of the model space}
By Theorem 2.1, the surfaces whose ends are asymptotically equal to 
$S \times(0,1)$ or $S \times(1,\infty),\, S:=S^1$,
equiped with the metric given by
\begin{equation}
ds^2 = \frac{(dy)^2 + (dx)^2}{y^2}, \quad 0 \leq x \leq 2 \pi,
\label{S3Modelmetric}
\end{equation}
form a broad and meaningfull class of 2-dimensional surfaces.
In this subsection, we shall introduce a  model for such  surfaces and study the
spectral properties of the Laplace-Beltrami operator on it. Since it is an unperturbed 
(free) space, we put the subscript {\it free} for every related object on it.
We put $M_{free} =S$ and let $\partial_x^2$ be the Laplace-Beltrami operator on $S$. It has eigenvalues and 
eigenvectors
\begin{equation}
\lambda_n =  n^2, \quad \varphi_n(x) = e^{ i nx}/\sqrt{2\pi }, \quad n \in {\bf Z}.\label{S3lambdanvarphin}
\end{equation}
Let $\mathcal M_{free} = M_{free}\times(0,\infty)$ and $H_{free}$  be given by  
\begin{equation}
 H_{free} = - y^2(\partial_y^2 + \partial_x^2)  - \frac{1}{4}.
\label{S3H0}
\end{equation}
$\mathcal M_{free}$ has two infinities corresponding to $y = 0$ and $y = \infty$. We call the former 
the {\it regular end},  and the latter the {\it cusp}. In the following, the subscripts $c$ and 
$reg$ mean the cusp and regular end, respectively.


\subsubsection{Green's operator}
 Green's kernel of $H_{free}$ is computed as follows.
Consider the 1-dimensional operators
\begin{equation}
L_{free}(\zeta) = y^2(- \partial_y^2 + \zeta^2) - \frac{1}{4}, \quad 
\zeta \in {\bf R},
\label{S3L0zeta}
\end{equation}
\begin{equation}
(L_{free}(\zeta) + \nu^2)^{-1} =: G_{free}(\zeta,\nu).
\label{S3L0zetaresolvent}
\end{equation}
If $\zeta \neq 0$, by (\ref{S3KnuandInu}), (\ref{S3BesselDiffEq}), $G_{free}(\zeta,\nu)$  
has the following expression (see \cite{IsKu10}, Chap. 1, \S 3),
\begin{equation}
\left(G_{free}(\zeta,\nu)\psi\right)(y) 
= \int_0^{\infty}G_{free}(y,y';\zeta,\nu)\psi(y')\frac{dy'}{(y')^2},
\label{S3G0GreenOp}
\end{equation}
\begin{equation}
G_{free}(y,y';\zeta,\nu) = \left\{
\begin{split}
\big(yy'\big)^{1/2}K_{\nu}(\zeta y)I_{\nu}(\zeta y'), \quad y > y' > 0, \\
\big(yy'\big)^{1/2}I_{\nu}(\zeta y)K_{\nu}(\zeta y'), \quad y' > y > 0.
\end{split}
\right.
\label{S3G0Kernel}
\end{equation}
Let us remark that in \cite{IsKu10}, $L_{free}$, $G_{free}$ are denoted by $L_0$, $G_0$. In what 
follows, the subscript 0 is, however, reserved to denote the terms associated with the 
eigenvalue $\lambda_0 = 0$.

When $\zeta = 0$, we have (see \cite{IsKu10}, Chap. 3, \S 2),
\begin{equation}
\left(G_{free}(0,\nu)\psi\right)(y) = \int_0^{\infty}G_{free}(y,y';0,\nu)\psi(y')\frac{dy'}{(y')^2},
\label{S30intop}
\end{equation}
\begin{equation}
G_{free}(y,y';0,\nu) = \frac{1}{2\nu}\left\{
\begin{split}
y^{\frac{1}{2}+\nu}(y')^{\frac{1}{2}-\nu}, \quad y' > y > 0, \\
y^{\frac{1}{2}-\nu}(y')^{\frac{1}{2}+\nu}, \quad y > y' > 0.
\end{split}
\right.
\label{S30kernel}
\end{equation}
We define $\mathcal B(\bf C)$ and $\mathcal B(\bf C)^{\ast}$ by putting ${\bf h} = {\bf C}$ in Subsection 3.1. Then we have, by \cite{IsKu10}, Chap. 1, Lemma 3.8,
\begin{equation}
\|G_{free}(\zeta,\nu)\psi\|_{\mathcal B(\bf C)^{\ast}} \leq C\|\psi\|_{\mathcal B(\bf C)},
\label{S3G0zetanuestimate}
\end{equation}
where the constant $C$ is independent of $\nu$, when $\nu$ varies over a compact set in $\{{\rm Re}\,\nu \geq 0\}\setminus{\bf Z}$, and also of $\zeta$, when ${\rm Re}\,\zeta > 0$. One can also prove (\ref{S3G0zetanuestimate}) for $\zeta = 0$.

Recalling (\ref{S3lambdanvarphin}), we put, for $f(x,y) \in 
{\mathcal H}_{free}:=L^2((0, \infty): L^2(S); dy/y^2)$,
\begin{equation}
\widehat f_n(y) = \int_{M_{free}}f(x,y)\overline{\varphi_n(x)}dx.
\label{S3hatfny}
\end{equation}
Let $R_{free}(z) = (H_{free} - z)^{-1}, \, z=-\nu^2$. Then
\begin{equation}
\begin{split}
 R_{free}(-\nu^2)f 
& = 
\sum_{n \in {\bf Z}}\varphi_n(x)
\left(\big(L_{free}(|n|) + \nu^2)^{-1}\widehat f_n(\cdot)\right)(y) \\
& = 
\sum_{n \in {\bf Z}}\varphi_n(x)
\left(G_{free}(|n|,\nu)\widehat f_n(\cdot)\right)(y).
\end{split}
\label{S3ResolventofH0}
\end{equation}
 For $0 < a < b$, we put
\begin{equation}
J_{\pm} = \{z \in {\bf C}\;;\; a \leq {\rm Re}\,z \leq b, \ \pm {\rm Im}\, z > 0\}.
\label{S3Jpm}
\end{equation}
The estimate (\ref{S3G0zetanuestimate}) then implies
\begin{equation}
\|R_{free}(z)f\|_{{\mathcal B}^{\ast}} \leq C\|f\|_{\mathcal B},
\label{S3R0zunifestimate}
\end{equation}
with $\mathcal B = \mathcal B(L^2(S))$ and $\mathcal B^{\ast} = \mathcal B(L^2(S))^{\ast}$, 
where the constant $C$ is independent of $z \in J_{\pm}$. 
This uniform estimate is crucial in proving the limiting absorption principle.
In fact, by \cite{IsKu10}, Chap. 3, Theorems 3.5 and 3.8, the following theorem holds. 


\begin{theorem} 
(1) $\sigma(H_{free}) = [0,\infty)$. \\
\noindent
(2) $\sigma_p(H_{free}) = \emptyset$. \\
\noindent
(3) For $\lambda > 0, f \in {\mathcal B} = \mathcal B( L^2(S))$, the following limit exists in the weak ${\ast}$-sense 
$$
\lim_{\epsilon \to 0}R_{free}(\lambda \pm i\epsilon)f =: 
R_{free}(\lambda \pm i0)f,
$$
i.e. there exits the limit
$$
\lim_{\epsilon\to 0}\left(R_{free}(\lambda \pm i\epsilon)f,g\right), \quad \forall f, g \in \mathcal B.
$$
\end{theorem}

Note that, since ${\rm Re}\, \nu \geq 0$, we have, letting $\nu = - i(k \pm i\epsilon)$, $k > 0$,
\begin{equation}
R_{free}(k^2 \pm i0)f = \sum_{n \in {\bf Z}}\varphi_n(x)
\left(G_{free}(|n|,\mp ik)\widehat f_n(\cdot)\right)(y).
\label{S3Rfreek2+i0}
\end{equation}


\subsubsection{Fourier transform}
Let $f \in C_0^{\infty}(\mathcal M_{free})$, and $k > 0$. 
For $n \neq 0$, the associated Fourier-Bessel transform  is defined by
\begin{equation}
{ F}_{free,n}(k)f = 
\frac{\big(2k\sinh(k\pi)\big)^{1/2}}{\pi}
\int_0^{\infty}y^{1/2}K_{ik}(|n| y)
\widehat f_n(y)\frac{dy}{y^2}.
\label{S3Fouriernneq0}
\end{equation}
For $n = 0$, the associated Mellin transform is defined by
\begin{equation}
F_{free,0}^{(\pm)}(k)f = \frac{1}{\sqrt{2\pi}}
 \int_0^{\infty}y^{\frac{1}{2} \pm ik}\widehat f_0(y)\frac{dy}{y^2}.
\label{Fouriern=0}
\end{equation}


\begin{definition} We put
\begin{equation}
{\bf h} = {\bf C}\oplus L^2(S),\quad 
\widehat{\mathcal H} = L^2((0,\infty) ; {\bf h} ; dk),
\nonumber
\end{equation} 
and define $\mathcal F^{(\pm)}_{c,free}(k)$ and $\mathcal F^{(\pm)}_{reg,free}(k)$ by
\begin{equation}
{\mathcal F}_{c,free}^{(\pm)}(k)f = 
F_{free,0}^{(\mp)}(k)f.
\label{S3Focplusminus}
\end{equation}
\begin{equation}
\begin{split}
\left({\mathcal F}_{reg,free}^{(\pm)}(k)f\right)(x)
   & =  C_0^{(\pm}(k)F_{free,0}^{(\pm)}(k)f \\
& + \sum_{n \in {\bf Z}\setminus\{0\}}
  C_{n}^{(\pm)}(k)\varphi_n(x)
F_{free,n}(k)f,
\end{split}
\label{S3F0regpm}
\end{equation}
\begin{equation}
C_{n}^{(\pm)}(k) = 
\left\{
\begin{split}
& \left(\frac{n}{2}\right)^{\mp ik} \quad (n \neq 0), \\
& \dfrac{\pm i}{k\omega_{\pm}(k)}\sqrt{\dfrac{\pi}{2}} \quad (n = 0),
\end{split}
\right.
\label{S3constCgammaastk}
\end{equation}
\begin{equation}
\omega_{\pm}(k) = \frac{\pi}{(2k\sinh(k\pi))^{1/2}\Gamma(1 \mp ik)}.
\label{S3omegak}
\end{equation}
Finally, we define the Fourier transform assocaited with $H_{free}$ by
\begin{equation}
{\mathcal F}_{free}^{(\pm)}(k)
 = \left({\mathcal F}_{c,free}^{(\pm)}(k),{\mathcal F}_{reg,free}^{(\pm)}(k)\right).
 \nonumber
\end{equation}
\end{definition}

The important step for the spectral representation is the following Parseval's formula
\begin{equation}
\frac{k}{\pi i}\left([R_{free}(k^2 + i0) - R_{free}(k^2 -i0)]f,f\right) = 
\|\mathcal F_{free}^{(\pm)}(k)f\|_{\bf h}^2.
\label{S3Parseval}
\end{equation}
This and the uniform estimate (\ref{S3R0zunifestimate}) imply the following inequality
\begin{equation}
\|\mathcal F_{free}^{(\pm)}(k)f\|_{\bf h} \leq C\|f\|_{\mathcal B}.
\label{S3Ffree(k)festimate}
\end{equation}
Therefore, $\mathcal F_{free}^{(\pm)}(k)$ can be extended uniquely on $\mathcal B$. For $f \in \mathcal B$, we define an ${\bf h}$-valued function of $k \in (0,\infty)$ by
$$
\left(\mathcal F_{free}^{(\pm)}f\right)(k) = \mathcal F_{free}^{(\pm)}(k)f.
$$
 Then, by integrating (\ref{S3Parseval}) with respect to $k$ over $(0,\infty)$, we see that $\mathcal F_{free}^{(\pm)}$ can be extended to an isometry from $\mathcal H_{free}$ to $\widehat{\mathcal H}$. In fact, it is unitary (see \cite{IsKu10}, Chap. 3, Theorem 2.5).


\begin{theorem}
 ${\mathcal F}_{free}^{(\pm)}$ is uniquely extended to a unitary operator from $\mathcal H_{free}$ to
$\widehat{\mathcal H}$. Moreover, if $f \in D(H_{free})$
\begin{equation}
({\mathcal F}_{free}^{(\pm)}H_{free}f)(k) = k^2
({\mathcal F}_{free}^{(\pm)}f)(k).
\nonumber
\end{equation}
\end{theorem}

\medskip
The Fourier transform $\mathcal F_{free}^{(\pm)}$ is related to the asymptotic expansion of the resolvent at infinity in the following way.


\begin{theorem}
For $k > 0$ and $f \in {\mathcal B}$, we have
\begin{equation}
\lim_{R\to\infty}\frac{1}{\log R}\int_{1/R}^1
\|\big(R_{free}(k^2 \pm i0)f\big)(\cdot,y) - v^{(\pm)}_{reg}(\cdot,y)\|_{L^2(S)}^2\frac{dy}{y^2} = 0,
\label{S3R0nearregular}
\end{equation}
\begin{equation}
v^{(\pm)}_{reg}(x,y) = \omega_{\pm}(k)\, y^{\frac{1}{2}  \mp ik}
\left({\mathcal F}^{(\pm)}_{reg,free}(k)f\right)(x),
\nonumber
\end{equation}
\begin{equation}
\lim_{R\to\infty}\frac{1}{\log R}\int_{1}^R
\|\big(R_{free}(k^2 \pm i0)f\big)(\cdot,y)
 - v_{c}^{(\pm)}\|_{L^2(S)}^2\frac{dy}{y^2} = 0,
\label{S3R0nearcusp}
\end{equation}
\begin{equation}
v_{c}^{(\pm)} = \omega_{\pm}^{(c)}(k)\,
y^{\frac{1}{2} \pm ik}
{\mathcal F}_{c,free}^{(\pm)}(k)f,
\nonumber
\end{equation}
where 
\begin{equation}
\omega_{\pm}^{(c)}(k) = \pm \frac{i}{k}\sqrt{\frac{\pi}{2}}.
\label{S3omegack}
\end{equation}
\end{theorem}

This theorem is proven by comparing the form of Green's function
(\ref{S3G0Kernel}), (\ref{S30kernel}) with the definition of $\mathcal F_{free}^{(\pm)}$, and using the asymptotic expansion of Bessel functions. See \cite{IsKu10}, Chap. 3, Theorem 2.6.


\subsection{Basic spectral properties for asymptotically hyperbolic manifolds}
We turn to the spectral properties of the manifold $\mathcal M$ satisfying the 
assumptions {\bf(A-1)} $\sim$ {\bf(A-4)} in \S 1. To deal with the Laplace-Beltrami 
operator $-\Delta_g$ for $\mathcal M$, we first pass it to the gauge transformation
\beq \label{2.3.1}
- \Delta_g - \frac{1}{4} \to - \rho^{1/4}\Delta_g\rho^{-1/4} - \frac{1}{4}.
\nonumber
\eeq
Here $\rho \in C^{\infty}(\mathcal M)$ is a positive function such that $\rho =1$ in a small neighborhood of $\mathcal M_{sing}$.
On each end $\mathcal M_j$,
\begin{equation}
\rho =  g_{free(j)}/ g,
\nonumber
\end{equation}
where $g_{free(j)}$ and $g$ define the volume elements, in the $(x, y)-$coordinates,
of the unperturbed and perturbed metrics on $\mathcal M_j$. Note that
$\rho=\hbox{const}$ in each ${\mathcal M}_i,\, i=1, \dots, \mu$. 
Let $H$ be the self-adjoint extension of $- \rho^{1/4}\Delta_g\rho^{-1/4} - 1/4$ defined 
in the same way as in Lemma 2.8. 
Our first concern is the (non) existence of 
the embedded eigenvalues in the continuous spectrum.


\begin{theorem} \label{th.3.6}
(1) $\sigma_{e}(H) = [0,\infty)$. \\
(2) If one of ${\mathcal M}_i$'s is a regular end, then $\sigma_p(H)\cap(0,\infty) = 
\emptyset$. \\
\noindent
(3) If all of the ${\mathcal M}_i$'s have a cusp, then $\sigma_p(H)\cap(0,\infty)$ is discrete with finite multiplicities, whose possible accumulation points are $0$ and $\infty$.
\end{theorem}

For the proof, see \cite{IsKu10}, Chap. 3, Theorems 3.2 and 3.5. The assertion (1) is a consquence of Theorem 3.2 (1) and Weyl's theorem on the perturbation of essential spectrum. The main tool for proving the assertion (2)  is  a theorem on the growth property of solutions to an abstract differential equation with operator-valued coefficients (\cite{IsKu10}, Chap. 2, Theorem 3.1). The assertion (3) is a standard result which follows from the a-priori estimates for solutions to the reduced wave equation 
\begin{equation}
\big(- \Delta_g - \frac{1}{4} - z\big)u = f
\label{S3Reduceswave}
\end{equation}
 and the short-range perturbation theory for the Schr{\"o}dinger equation. \qed

\medskip
Take $\chi_{0} \in C_0^{\infty}(\mathcal M)$ such that $\chi_{0} = 1$ on $\mathcal K$, and put $\chi_i = 1 - \chi_{0}$ on $\mathcal M_i$, $\chi_i = 0$ on $\mathcal M \setminus \mathcal M_i$. Then $\{\chi_{0}, \chi_1,\cdots,\chi_N\}$ is a partition of unity on $\mathcal M$
{\nntext subordinated to decomposition (\ref{1.0}).}

We define the Besov space $\mathcal B_i$ by $\mathcal B_i = \mathcal B(\bf C)$, when $\mathcal M_i$ has a cusp, and $\mathcal B_i = 
\mathcal B( L^2(S^{r_i}))$, when $\mathcal M_i$ has a regular infinity.  We then put
$$
\|f\|_{\mathcal B} = \|\chi_{0}f\|_{L^2(\mathcal M)} + \sum_{i=1}^N\|\chi_if\|_{\mathcal B_i},
$$
$$
\|u\|_{\mathcal B^{\ast}} = \|\chi_{0}u\|_{L^2(\mathcal M)} + \sum_{i=1}^N\|\chi_iu\|_{{\mathcal B_i}^{\ast}},
$$
which define the Besov type spaces $\mathcal B$ and $\mathcal B^{\ast}$ on $\mathcal M$. 

Let $R(z) = (H - z)^{-1}$ be the resolvent of $H$.

 
\begin{theorem}
For $\lambda \in (0,\infty)\setminus\sigma_p(H)$, there exists a limit
$$
\lim_{\epsilon \to 0}R(\lambda \pm i\epsilon) \equiv 
R(\lambda \pm i0) \in {\bf B}({\mathcal B}\,;{\mathcal B}^{\ast})
$$
in the weak $\ast$-sense. Moreover, for any compact interval $I \subset (0,\infty)\setminus\sigma_p(H)$, there exists a constant $C > 0$ such that
\begin{equation}
\|R(\lambda \pm i0)f\|_{{\mathcal B}^{\ast}} \leq C\|f\|_{\mathcal B},
\quad \lambda \in I.
\nonumber
\end{equation} 
 For $f, g \in {\mathcal B}$,  $(R(\lambda \pm i0)f,g)$ is continuous with respect to $\lambda \in (0,\infty)\setminus\sigma_p(H)$.
\end{theorem}

This theorem is proved in \cite{IsKu10}, Chap. 3, Theorem 3.8. The proof consists of 
two 
main ingredients. We first establish some a-prori estimates for solutions to the reduced 
wave equation  (\ref{S3Reduceswave}) by the elementary tool of integration by parts 
(\cite{IsKu10}, Chap. 2, Lemmas 2.4 $\sim$ 2.8). This 1st step is essentially the 
1-dimensional problem. 
The proof of Theorem 3.7 is done by the argument of contradiction, using
 the compactness of the perturbation and reducing the problem to the uniqueness 
 of solutions of the equation (\ref{S3Reduceswave}) satisfying the 
 corresponding radiation condition.

The above mentioned radiation condition is as follows. Let 
$$
\sigma_{\pm}(\lambda) = \frac{1}{2} \mp i\sqrt{\lambda}, \quad \lambda > 0.
$$
We say that a solution $u \in \mathcal B^{\ast}$ of the equation $(- \Delta_g - \frac{1}{4} - \lambda)u = f \in \mathcal B$ satisfies the {\it outgoing radiation condition}, or $u$ is {\it outgoing}, if 
\begin{equation}
\begin{split}
\label{S3RadCond}
&\displaystyle \frac{1}{\log R}\int_{2}^{R}\|\big(y\partial_y - \sigma_+(\lambda)\big)u(\cdot,y)\|^2_{L^2(S^{r_j})}\frac{dy}{y^2} \to 0, 
\quad (j = 1, \cdots, \mu), \\
& \displaystyle \frac{1}{\log R}\int_{1/R}^{1/2}\|\big(y\partial_y - \sigma_+(\lambda)\big)u(\cdot,y)\|^2_{L^2(S^{r_j})}\frac{dy}{y^2} \to 0,
\quad (j = \mu+1, \cdots, N)
\end{split}
\end{equation}
hold as $R \to \infty$.
The following theorem follows from \cite{IsKu10}, Chap. 3, Theorems 3.7 and 3.8.

\begin{theorem}
Let $\lambda \in (0,\infty)\setminus\sigma_p(H)$. \\
\noindent
(1) If $u \in \mathcal B^{\ast}$ satisfies $(H - \lambda) u = 0$ and is outgoing, then $u = 0$. \\
\noindent
(2) For $f \in \mathcal B$, $R(\lambda + i0)f$ is outgoing.
\end{theorem}


\subsection{Fourier transforms associated with $H$} 
We shall make use of the perturbation method to construct the Fourier transform for $H$ from that of the model space.
Let $H_{free(j)}$ be defined by
\begin{equation}
H_{free(j)} = - y^2(\partial_y^2 + \Delta_{M_j}) - \frac{1}{4},
\label{S3Hfreej}
\end{equation}
where $\Delta_{M_j}$ is the Laplace-Beltrami operator of $M_j$. Let
$\chi_j$ be the partition of unity as above. We put
\begin{equation}
\widetilde V_j = H - H_{free(j)} \quad {\rm on} \quad \mathcal M_j.
\label{C3S3widetildeVj}
\end{equation}
This is symmetric on $C_0^{\infty}(\mathcal M_j)$, since so are $H$ and $H_{free(j)}$.
Observe that
\begin{equation}
(H_{free(j)} - \lambda)\chi_j 
Q_j(\la \pm i0)
R(\lambda \pm i0),
\nonumber
\end{equation}
where
\begin{equation}
Q_j(z) = \chi_j + \left([H_{free(j)},\chi_j] - \chi_j\widetilde V_j\right)R(z).
\label{C3S3Qjz}
\end{equation}
Therefore, we have the following equality
\begin{equation}
 \chi_jR(\lambda \pm i0) 
 =  
R_{free(j)}(\lambda \pm i0) Q_j(\lambda \pm i0). 
\label{S3RandRj}
\end{equation}
This formula suggests how the generalized Fourier transform is constructed by 
the perturbation method.

Let $\lambda_{j,n} = (n/r_j)^2$, $\varphi_{j,n}(x) = e^{inx/r_j}/\sqrt{2\pi r_j}$ be the eigenvalues and normalized eigenvectors of $\Delta_{M_j}$. 
We  define $\mathcal F_{c,free(j)}^{(\pm)}(k)$ by (\ref{S3Focplusminus}), 
and  $\mathcal F^{(\pm)}_{reg,free(j)}(k)$  by (\ref{S3F0regpm}) with $M$ 
replaced 
by $M_j$, $\varphi_{n}$ by 
$\varphi_{j,n}$, and 
$C_n^{(\pm)}(k)$ by $C_{j,n}^{(\pm)}(k)$, i.e.
\begin{equation}
C_{j,n}^{(\pm)}(k) = 
\left\{
\begin{split}
& \left(\frac{\sqrt{\lambda_{j,n}}}{2}\right)^{\mp ik}, \quad
(\lambda_{j,n} \neq 0), \\
& \frac{\pm i}{k\omega_{\pm}(k)}\sqrt{\frac{\pi}{2}}, \quad (\lambda_{j,n} = 0).
\end{split}
\right.
\label{C3S3Cjmk}
\end{equation}

\subsubsection{Definition of ${\mathcal F}_{free(j)}^{(\pm)}(k)$} 
Recall that, for $1 \leq j \leq \mu$, $\mathcal M_j$ has a cusp, and, for $\mu + 1 \leq j \leq N$, $\mathcal M_j$ has a regular infinity.

\medskip
\noindent
(i) For $1 \leq j \leq \mu$ (the case of cusp), we define
\begin{equation}
{\mathcal F}_{free(j)}^{(\pm)}(k) = \mathcal F_{c,free(j)}^{(\pm)}(k).
\label{C3S3Ffreejpmkcusp} 
\end{equation}

\medskip
\noindent 
(ii) For $\mu + 1 \leq j \leq N$ (the case of regular infinity), we define
\begin{equation}
{\mathcal F}_{free(j)}^{(\pm)}(k) = 
\mathcal F^{(\pm)}_{reg,free(j)}(k),
\label{C3S3Ffreejpmk}
\end{equation}

\subsubsection{Definition of ${\mathcal F}^{(\pm)}(k)$}
 For $1 \leq j \leq N$, we define
\begin{equation}
\mathcal F^{(\pm)}_j(k) = \mathcal F^{(\pm)}_{free(j)}(k)Q_j(k^2 \pm i0),
\label{C3S3Fpmjkdefine}
\end{equation}
Finally, we define the Fourier transform associated with $H$ by
\begin{equation}
{\mathcal F}^{(\pm)}(k) = \big({\mathcal F}_{1}^{(\pm)}(k), \cdots, {\mathcal F}_{N}^{(\pm)}(k)\big).
\label{C3S3Fpmkdefine}
\end{equation}

\subsubsection{Eigenfunction expansion theorem}
Let
\begin{equation}
{\bf h}_{\infty} = \oplus_{j= 1}^{N} {\bf h}_j, \quad {\bf h}_j={\bf C},\, 1 \leq j \leq \mu,
\quad {\bf h}_j=L^2(M_j), \,\mu+ 1 \leq j \leq N.
\label{eq:Chap3Sect2hinfty}
\end{equation}
and for $\varphi, \psi \in {\bf h}_{\infty}$, define the inner product by
\begin{equation}
(\varphi,\psi)_{{\bf h}_{\infty}} = \sum_{j=1}^{\mu}\varphi_j\overline{\psi_j}|M_j| + \sum_{j=\mu+1}^{N}(\varphi_j,\psi_j)_{L^2(M_j)},
\label{C3S3hinftyinnerproduct}
\end{equation}
where $|M_j| = 2\pi r_j$ is the length of $M_j$.
We put
\begin{equation}
\widehat{\mathcal H} = L^2((0,\infty);{\bf h}_{\infty};dk).
\nonumber
\end{equation}


\begin{theorem}
We define $\big({\mathcal F}^{(\pm)}f\big)(k) = {\mathcal F}^{(\pm)}(k)f$ for $f \in {\mathcal B}$. Then  ${\mathcal F}^{(\pm)}$ is uniquely extended to a bounded operator from $L^2({\mathcal M})$ to $\widehat{\mathcal H}$ with the following properties.

\noindent
(1) $\ {\rm Ran}\,\,{\mathcal F}^{(\pm)} = \widehat{\mathcal H}$. \\
\noindent
(2) $\ \|f\| = \|{\mathcal F}^{(\pm)}f\|$ for $f \in {\mathcal H}_{ac}(H)$. \\
\noindent
(3) $\ {\mathcal F}^{(\pm)}f = 0$ for $f \in {\mathcal H}_p(H)$.  \\
\noindent
(4) $\ 
\left({\mathcal F}^{(\pm)}Hf\right)(k) = 
k^2\left({\mathcal F}^{(\pm)}f\right)(k)$ for 
$f \in D(H)$. \\
\noindent
(5)$\ {\mathcal F}^{(\pm)}(k)^{\ast} \in {\bf B}({\bf h}_{\infty};{\mathcal B}^{\ast})$ and
$(H - k^2){\mathcal F}^{(\pm)}(k)^{\ast} = 0$  for $k^2 \in (0,\infty)\setminus\sigma_p(H)$. \\
\noindent
(6) For $f \in {\mathcal H}_{ac}(H)$, the inversion formula holds:
\begin{eqnarray*}
f = \left({\mathcal F}^{(\pm)}\right)^{\ast}{\mathcal F}^{(\pm)}f 
= \sum_{j=1}^N\int_0^{\infty}{\mathcal F}^{(\pm)}_j(k)^{\ast}
\left({\mathcal F}^{(\pm)}_jf\right)(k)dk.
\end{eqnarray*}
\end{theorem}

The most important step of the proof of this theorem is  Parseval's formula
\begin{equation}
\frac{k}{\pi i}\left(\big[R(k^2 + i0) - R(k^2 - i0)\big]f,g\right) = 
\left(\mathcal F^{(\pm)}(k)f,\mathcal F^{(\pm)}(k)g\right)_{{\bf h}_{\infty}}
\nonumber
\end{equation}
for $f, g \in \mathcal B$, $k^2 \in (0,\infty)\setminus\sigma_p(H)$ (\cite{IsKu10}, Chap. 3, Lemma 3.11), which is proven by the following Theorem 3.10. The remaining arguments are routine.
See \cite{IsKu10}, Chap. 3, Theorem 3.12 for the details.

\noindent
{\it Remark 1}. $\ $ The meaning of the integral in (6) is as follows. Let $(0,\infty)\setminus\sigma_p(H) = \cup_{i=1}^{\infty}I_i$, where $I_i = (a_i,b_i)$ 
are non-overlapping open intervals. For $g(k) \in \widehat{\mathcal H}$, we have by (5)
$$
\int_{\sqrt{a_i}+\epsilon}^{\sqrt{b_i}-\epsilon}{\mathcal F}^{(\pm)}_j(k)^{\ast}
g(k)dk \in {\mathcal B}^{\ast}.
$$
As a matter of fact, it belongs to $L^2({\mathcal M})$, and 
$$
\lim_{\epsilon \to 0}\int_{\sqrt{a_i}+\epsilon}^{\sqrt{b_i}-\epsilon}{\mathcal F}^{(\pm)}_j(k)^{\ast}
g(k)dk \in L^2({\mathcal M})
$$
in the sense of strong convergence in $L^2({\mathcal M})$. Denoting this limit by
$$
\int_{\sqrt{I_i}}{\mathcal F}^{(\pm)}_j(k)^{\ast}
g(k)dk,
$$
we define
$$
\int_{0}^{\infty}{\mathcal F}^{(\pm)}_j(k)^{\ast}
g(k)dk = 
\sum_{i=1}^{\infty}\int_{\sqrt{I_i}}{\mathcal F}^{(\pm)}_j(k)^{\ast}
g(k)dk.
$$

\subsubsection{Asymptotic expansion of the resolvent}
For $f, g \in {\mathcal B}^{\ast}$ on ${\mathcal M}$, by
$f \simeq g$
we mean that on each end the following relation holds,
\begin{equation}
\lim_{R\to\infty}\frac{1}{\log R}\int_{1/R}^{R}\rho_j(y)\|f(y) - g(y)\|^2_{L^2(M_j)}
\frac{dy}{y^2} = 0
\nonumber
\end{equation}
 where $\rho_j(y) = 1$ $(y<1/2)$, $\rho_j(y) = 0$ $(y > 1)$, when $\mathcal M_j$ has a 
regular infinity, and $\rho_j(y) = 0$ $(y < 1)$, $\rho_j(y) = 1$ $(y > 2)$, when $\mathcal M_j$ 
has a cusp. 
Theorem 3.5 shows that $\mathcal F_{free(j)}^{(\pm)}(k)f$ is computed from the asymptotic 
expanison of $R_{free(j)}(\lambda\pm i0)f$ at infinity. This, combined with the formula (\ref{S3RandRj}) and  definition  (\ref{C3S3Fpmjkdefine}), implies the following 
theorem (see \cite{IsKu10}, Chap. 3, Theorem 3.10).


\begin{theorem} Let $f \in {\mathcal B}$, 
$k^2 \in \sigma_e(H)\setminus\sigma_p(H)$, and $\chi_j$ be the partition of unity on $\mathcal M$. 
Then we have
\begin{equation}
\begin{split}
R(k^2 \pm i0) f &\simeq 
 \omega_{\pm}^{(c)}(k)\sum_{j=1}^{\mu}
\chi_jy^{1/2\pm ik}
{\mathcal F}^{(\pm)}_j(k)f \\
& + \omega_{\pm}(k)\sum_{j=\mu + 1}^N
\chi_jy^{1/2\mp ik}
{\mathcal F}^{(\pm)}_j(k)f.
\end{split}
\nonumber
\end{equation}
\end{theorem}

The following theorem is a characterization of the solution space of the Helmholtz equation, and is proved in the same way as in \cite{IsKu10}, Chap. 2, Theorem 7.8.


\begin{theorem}
If $k^2 \in (0,\infty)\setminus\sigma_p(H)$, we have 
\begin{equation}
{\mathcal F}^{(\pm)}(k){\mathcal B} = {\bf h}_{\infty},
\nonumber
\end{equation}
\begin{equation}
\{u \in {\mathcal B}^{\ast}\,;\,(H - k^2)u = 0\} = 
{\mathcal F}^{(\pm)}(k)^{\ast}{\bf h}_{\infty}.
\nonumber
\end{equation}
\end{theorem}


\subsection{$S$ matrix} 

We derive an asymptotic expansion of solutions to the Helmholtz equation. Let $V_{\ell}$ be the differential operator defined by
\begin{equation}
 V_\ell =  [H_{free(\ell)},\chi_\ell] - \chi_\ell \widetilde V_\ell \quad 
 (1 \leq \ell \leq N),
 \nonumber
\end{equation}
where $\widetilde V_\ell$ is defined by (\ref{C3S3widetildeVj}). 
We put
\begin{equation}
J_j(k) = \sum_{\la_{j, m} \neq 0}
\left(\frac{\sqrt{\lambda_{j,m}}}{2}\right)^{-2ik}
P_{j,m} = 
\left(\frac{\sqrt{-\Delta_{M_j}}}{2}\right)^{-2ik}P_j^+,
\label{C3S3Jpk}
\end{equation}
where $\Delta_{M_j}$ is the Laplace-Beltami operator on $M_j$ and 
$P_j^+$ is the projection onto the subspace on which $- \Delta_{M_j} > 0$.
For $1 \leq j, \ell \leq N$, we define 
$\widehat S_{j \ell}(k) \in {\bf h}_{\ell};\,{\bf h}_{j}$ by
\begin{equation}
 \widehat S_{j\ell }(k) = 
 \left\{
 \begin{split}
 &  \frac{\pi i}{k}
 {\mathcal F}_j^{(+)}(k)\big(V_\ell\big)^{\ast}\left({\mathcal F}_{free(\ell)}^{(-)}(k)\right)^{\ast},  \quad 1 \leq j \leq \mu,
\\
& \delta_{j \ell}J_j(k) + \frac{\pi i}{k}{\mathcal F}_j^{(+)}(k)\big(V_\ell\big)^{\ast}\left({\mathcal F}_{free(\ell)}^{(-)}(k)\right)^{\ast}, 
\quad \mu + 1 \leq j \leq N.
\end{split}
\right.
\label{C3S3widehatSpqk}
\end{equation}

We define an operator-valued $N \times N$ matrix $\widehat S(k)$ by
\begin{equation}
\widehat S(k) = \Big(\widehat S_{j \ell}(k)\Big)_{j, \ell=1}^N,
\label{C3S3Smatrix}
\end{equation}
and call it  {\it $S$-matrix}. This is a bounded operator on ${\bf h}_{\infty}$.


\begin{theorem}
(1)  For any $u \in {\mathcal B}^{\ast}$ satisfying $(H - k^2)u = 0$, there exists a unique 
$\psi^{(\pm)} = (\psi^{(\pm)}_1,\cdots,\psi^{(\pm)}_N) \in {\bf h}_{\infty}$ such that
\begin{equation}
 \begin{split}
    u &\simeq \omega_-^{(c)}(k)\sum_{j=1}^{\mu}
   \,\chi_j \, y^{1/2-ik}\, \psi_{j}^{(-)}  + \omega_-(k)\sum_{j=\mu+1}^N\,\chi_j\, 
  y^{1/2+ik}\,\psi_j^{(-)} \\
   & -  \omega_+^{(c)}(k)\sum_{j=1}^{\mu}
  \,\chi_j\, y^{1/2+ik}\,
  \psi_{j}^{(+)} -  \omega_+(k)\sum_{j=\mu+1}^N
  \chi_j\, y^{1/2-ik}\,
  \psi_j^{(+)} .
   \end{split}
   \nonumber
\end{equation}
(2) For any $\psi^{(-)} \in {\bf h}_{\infty}$, there exists a unique $\psi^{(+)} \in {\bf h}_{\infty}$ and $u \in \mathcal B^{\ast}$ satisfying $(H - k^2)u = 0$, for which the expansion (1) holds. Moreover
\begin{equation}
\psi^{(+)} = \widehat S(k)\psi^{(-)}.
\nonumber
\end{equation}
(3)
$\widehat S(k)$ is unitary on ${\bf h}_{\infty}$.
\end{theorem}

For the proof, see \cite{IsKu10}, Chap. 3, Theorems 3.14, 3.15, 3.16.


\subsection{Helgason's theorem}
Before closing this section, we give some remarks on Theorems 3.11 and 3.12. 
As the most fundamental example of hyperbolic space, let us consider the Poincar{\'e} disc $D$ in ${\bf C}$. As is well-known, the Poisson integral
\begin{equation}
u(z) = \frac{1}{2\pi}\int_0^{2\pi}
\left(\frac{1 - |z|^2}{|e^{i\theta}-z|^2}\right)^sf(\theta)d\theta,
\label{S3Poisson}
\end{equation}
$f(\theta)$ being a function on the boundary $\partial D = S^1$,
gives a solution to the Helmholtz equation in $D$:
\begin{equation}
(- \Delta_g - E)u = 0,  \quad
E = 4s(s-1).
\label{S3EqPoincare}
\end{equation}
Our solution space $\mathcal B^{\ast}$, which is associated with the case in which the 
boundary space is $L^2(S^1)$, 
 has the following feature:
Regarding the decay at infinity, which corresponds to the boundary $\p D =  S^1$, 
of solutions for (\ref{S3EqPoincare}), $\mathcal B^{\ast}$ 
is the smallest space. In fact, by \cite{IsKu10}, Chap. 3, Theorem 3.6, if a solution 
$u$ of the equation (\ref{S3EqPoincare}) has a faster decay rate than $\mathcal B^{\ast}$ 
at regular infinity, $u$ vanishes identically. The largest solution space for 
(\ref{S3EqPoincare}) was given by Helgason. 
In \cite {Hel70}, he proved that {\it all} solutions of the Helmholtz equation
 is written by (\ref{S3Poisson}), where $f(\theta)$ is Sato's hyperfunction on the boundary. This result was extended to real hyperbolic spaces by \cite{Mine75} and to general symmetric spaces of rank 1 by \cite{KKMOOT78}.

\medskip
\noindent
{\it Remark 2.} Let $A(S^1)$ be the space of functions on $S^1$ having analytic continuations in a neighborhood of $S^1$. By the correspondence
\begin{equation}
{\bf c} = \left(c_n\right)_{n\in{\bf Z}} \Longleftrightarrow f_{\bf c} = 
\sum_{n\in{\bf Z}}c_ne^{ inx},
\label{S3candFourier}
\end{equation}
$A(S^1)$ is identified with the set of sequences
\begin{equation}
{\bf c}:\quad
\exists \rho > 1 \quad {\rm s.t.} \quad \sum_{n\in{\bf Z}}|c_n|\rho^{|n|} < \infty.
\nonumber
\end{equation}
The dual space of $A(S^1)$, the space of Sato's hyperfunctions on $S^1$, is identified with the set of sequences
\begin{equation}
{\bf d}=\left(d_n\right)_{n\in{\bf Z}}:\quad
 0< \forall\rho < 1, \quad \sup_{n\in{\bf Z}}|d_n|\rho^{|n|} < \infty.
\nonumber
\end{equation}

 Although $\mathcal B^{\ast}$ is the smallest solution space, it has sufficiently many 
 solutions if 
 one of the ends is regular. In fact, one can determine the whole manifold from the knowledge of 
 a component of the S-matrix associated with regular end, see \cite{IsKu10}. It is not the case for the cusp due to 
 the fact that the cusp gives rise only to  the 1-dimensional contribution to
  the continuous spectrum. 
 This requires us to generalize the notion of the S-matrix.


\section{Generalized S-matrix}


\subsection{Exponentially growing solutions}
In order to enlarge the solution space of the Helmholtz equation, we enlarge the associated space at infinity.

\begin{definition}
We introduce the sequential spaces $l^{2,\pm \infty}$ by
\begin{equation}
l^{2,\infty} \ni {\bf a} = (a_n)_{n\in{\bf Z}} \Longleftrightarrow 
 \forall \rho > 1, \quad
\sum_{n\in{\bf Z}}|a_n|^2\rho^{|n|} < \infty,
\nonumber
\end{equation}
\begin{equation}
l^{2,-\infty} \ni {\bf b} = (b_n)_{n\in{\bf Z}} \Longleftrightarrow 
\exists \rho > 1, \quad
\sum_{n\in{\bf Z}}|b_n|^2\rho^{-|n|} < \infty.
\nonumber
\end{equation}
\end{definition}

By the correspondence (\ref{S3candFourier}),  
$l^{2,\infty}$ is identified with the space of functions on $S^1$ having analytic continuations on ${\bf C}\setminus\{0\}$, moreover
\begin{equation}
l^{2,\infty} \subset A(S^1), \quad A(S^1)' \subset l^{2,-\infty}.
\nonumber
\end{equation}

Let $0 \neq k \in {\bf R}$. Suppose $u(x,y) \in C^{\infty}({\bf R} \times (1, \infty))$ is
$2 \pi r$-periodic in $x$, $u(x,y) = u(x+2 \pi r,y)$, and satisfies there 
the equation
\begin{equation}
- y^2\big(\partial_x^2 + \partial_y^2\big)u - \frac{1}{4}u = k^2 y.
\label{S4freeeq}
\end{equation}
 Expanding $u$ into the Fourier series
$$
u(x,y) = \frac{1}{\sqrt{2\pi r}}\sum_{n\in{\bf Z}}e^{ inx/r}u_n(y),
$$
we have
$$
y^2\Big(- \partial_y^2 + \frac{n^2}{r^2}\Big)u_n(y) - \frac{1}{4}u_n(y) = k^2 u_n(y), \quad y > 1.
$$
Then $u_n$ is written as
\begin{equation}
u_n(y) = \left\{
\begin{split}
& a_n\, y^{1/2}I_{-ik}(|n|y/r) + b_n\, y^{1/2}K_{ik}(|n|y/r), \quad (n\neq 0), \\
& a_0\, y^{1/2 - ik} + b_0\,y^{1/2 + ik}, \quad (n = 0).
\end{split}
\right.
\label{S4anbn}
\end{equation}
Let us note that  $K_{-\nu}(z) = K_{\nu}(z)$.


\begin{lemma}
Given $u(x,y) \in C^{\infty}({\bf R} \times (1, \infty))$, which is 
$2 \pi r-$periodic in $x$ 
 and satisfies (\ref{S4freeeq}), let ${\bf a} = (a_n)_{n\in{\bf Z}}$, ${\bf b} = (b_n)_{n\in{\bf Z}}$ be defined by (\ref{S4anbn}).
If ${\bf a} \in l^{2,\infty}$, then ${\bf b} \in l^{2,-\infty}$.
\end{lemma}
\noindent
{\bf Proof.} 
Recall the asymptotic expansion of the modified Bessel functions (\ref{S3Inunearinfty}), (\ref{S3Knunearinfty}).
Since ${\bf a} \in l^{2,\infty}$, we have $\sum_n|a_n|^2\big|I_{-ik}(|n|y/r)\big|^2 < \infty$ for any $y > 1$. By Parseval's formula,
\begin{equation}
\begin{split}
& y^{-1}\|u(\cdot,y)\|^2_{L^2(0,2\pi r)} \\
& = \sum_{n\neq 0}\big|a_nI_{-ik}(|n|y/r) + b_nK_{ik}(|n|y/r)\big|^2 + |a_0y^{-ik} + b_0y^{ik}|^2.
\end{split}
\nonumber
\end{equation}
We then have $\sum_{n\neq0}|b_n|^2\big|K_{ik}(|n|y/r)\big|^2 < \infty,\, y>1$, hence 
${\bf b} \in l^{2,-\infty}$. \qed

\medskip
We introduce the {\it spaces of generalized scattering data at infinity} :
\begin{equation}
{\bf A}_{\pm\infty} = \left({\mathop\oplus_{j=1}^{\mu}}l^{2,\pm\infty}\right)\oplus
\left({\mathop\oplus_{j=\mu+1}^N}L^2(M_j)\right),
\label{S4bfApm}
\end{equation}
$M_j$ being $S^{r_j}$ with metric $ds^2 = (dx)^2, \, x \in [0, 2 \pi r_j)$.

\medskip
We use the following notation. For
\begin{equation}
\psi^{(in)} = ({\bf a}_1,\cdots,{\bf a}_{\mu},\psi_{\mu+1}^{(in)},\cdots,\psi^{(in)}_N) \in {\bf A}_{\infty},
\label{S4psi-}
\end{equation}
\begin{equation}
\psi^{(out)} = ({\bf b_1},\cdots,{\bf b}_{\mu},\psi_{\mu+1}^{(out)},\cdots,\psi^{(out)}_N) \in {\bf A}_{-\infty},
\label{S4psi+}
\end{equation}
we let
\begin{equation}
u_j^{(in)} = 
\left\{
\begin{split}
&a_{j,0}\,y^{1/2-ik} + \sum_{n\neq0}
a_{j,n}\,e^{ inx/r_j}y^{1/2}I_{-ik}(|n|y/r_j), \quad 1 \leq j \leq \mu\\
&\omega_-(k)\,y^{1/2+ik}\psi_j^{(in)}(x), \quad \mu +1 \leq j \leq N,
\end{split}
\right.
\label{S4uj-}
\end{equation}
\begin{equation}
u_j^{(out)} = 
\left\{
\begin{split}
&b_{j,0}\,y^{1/2+ik} + \sum_{n\neq0}
b_{j,n}\,e^{ inx/r_j}y^{1/2}K_{ik}(|n|y/r_j), \quad 1 \leq j \leq \mu\\
&\omega_+(k)\, y^{1/2-ik}\psi_j^{(out)}(x), \quad \mu + 1 \leq j \leq N.
\end{split}
\right.
\label{S4uj+}
\end{equation}
Here $a_{j,n}, b_{j,n}$ are the $n$-th components of ${\bf a}_j \in l^{2,\infty}, {\bf b}_j \in l^{2,-\infty}$.
Let $\langle \; ,\; \rangle_j$ be the inner product of $L^2(S^{r_j})$ :
$$
\langle f,g\rangle_j = \int_{S^{r_j}}f\overline{g}\, dl.
$$


\begin{lemma}
Let $k > 0$ be such that $k^2 \not \in \sigma_p(H)$,  $\psi^{(in)}$, $u_j^{(in)}$ as in (\ref{S4psi-}), (\ref{S4uj-}), and
$u^{(in)} = \sum_{j=1}^N\chi_ju_j^{(in)}$. 
Then, there exists a unique solution $u$ such that 
\begin{equation}
(H - k^2) u = 0, \quad u - u^{(in)} \ {\rm is} \ {\rm outgoing}, 
\label{S4constructu}
\end{equation}
i.e. $u - u^{(in)}$ belongs to $\mathcal B^{\ast}$ on $\mathcal M$, and satisfies (\ref{S3RadCond}).
For this $u$, there exists $\psi^{(out)} = ({\bf b_1},\cdots,{\bf b}_{\mu},\psi_{\mu+1}^{(out)},\cdots,\psi^{(out)}_N) \in {\bf A}_{-\infty}$ such that \\
\noindent 
(1) for $j = 1,\cdots,{\mu}$,  
\begin{equation}
u = u_j^{(in)} - u_j^{(out)},  \quad \hbox{in} \quad \mathcal M_j\cap ({\rm supp}\,\chi_0)^c,
\label{S4uu-h+}
\end{equation}
(2) for $j = \mu+1,\cdots,N$,
\begin{equation}
u - u_j^{(in)} \simeq - u_j^{(out)}, \quad {\rm in} \quad \mathcal M_j.
\label{S4uminusu(-)}
\end{equation}
Explicitly, ${\bf b}_j$ and $\psi_j^{(out)}$ are given by
\begin{equation}
b_{j,0} = \frac{1}{2ik\sqrt{2\pi r_j}}\int_0^{\infty}(y)^{1/2-ik}\widehat f_{j,0}(y)\frac{dy}{(y)^2}, 
\label{S4bj0}
\end{equation}
\begin{equation}
b_{j,n} = \frac{1}{\sqrt{2\pi r_j}}\int_0^{\infty}y^{1/2}I_{-ik}(|n|y/d_j)\widehat f_{j,n}(y)\frac{dy}{y^2}, \quad n \neq 0,
\label{S4bjn}
\end{equation}

\begin{equation}
\psi_j^{(out)} = \mathcal F_j^{(+)}(k)f, \quad \mu + 1 \leq j \leq N,
\label{S4psijout}
\end{equation}
where
\begin{equation}
 f = (H-k^2)u^{(in)}, \quad
f_j = \chi_jf + [H_{free(j)},\chi_j]R(k^2 + i0)f,
\label{S4fandfj}
\end{equation}
\begin{equation}
\widehat f_{j,n} = \frac{1}{\sqrt{2\pi r_j}}\langle f_j,e^{ inx/r_j}\rangle_j.
\label{S4fjn}
\end{equation}
\end{lemma}

\noindent
{\bf Proof.} The uniqueness follows from Theorem 3.8. To show the existence,
we represent
\begin{equation}
u = u^{(in)} - R(k^2 + i0)f.
\label{S4u}
\end{equation}
Then the condition (\ref{S4constructu}) is satisfied by Theorem 3.8.
 By Theorem 3.10, we  have
\begin{equation}
R(k^2 + i0)f  \simeq  \ \omega_+^{(c)}(k)\sum_{j=1}^\mu\chi_jy^{1/2+ik}\mathcal F_j^{(+)}(k)f 
 +  \ \omega_+(k)\sum_{j=\mu+1}^N\chi_jy^{1/2-ik}\mathcal F_j^{(+)}(k)f,
\nonumber
\end{equation}
which proves (\ref{S4uminusu(-)}) and (\ref{S4psijout}).

For $j = 1,\cdots,\mu$, let
$H_{free(j)} = - y^2\Delta - 1/4$ on $S^{r_j}\times(0,\infty)$,
and put
\begin{equation}
R_{free(j)}(z) = (H_{free(j)} - z)^{-1}.
\nonumber
\end{equation}
Since
$$
(H_{free(j)} - \lambda)\chi_jR(\lambda \pm i0) = \chi_j + 
[H_{free(j)},\chi_j]R(\lambda \pm i0),
$$
we have 
\begin{equation}
\chi_jR(\lambda \pm i0) = R_{free(j)}(\lambda \pm i0)\chi_j  + 
R_{free(j)}(\lambda \pm i0)[H_{free(j)},\chi_j]R(\lambda \pm i0).
\label{S4chijresolvent}
\end{equation}
Note that on $\mathcal M_j$,
$f =  [H,\chi_j]u_j^{(-)}$,
and $[H,\chi_j]$ is a 1st-order differential operator with coefficients which are compactly 
supported in $\mathcal M_j$.
Therefore, $f_j$ is compactly supported, 
in particular $f_j = 0$ on $\mathcal M_j\cap\left({\rm supp}\,\chi_0\right)^c$, and, by 
(\ref{S4fandfj}),
\begin{equation}
\chi_jR(k^2 + i0)f = R_{free(j)}(k^2 + i0)f_j.
\label{S4chijandfj}
\end{equation}
By (\ref{S3G0GreenOp}), (\ref{S3G0Kernel}),
(\ref{S3ResolventofH0}), and taking account of (\ref{S3Rfreek2+i0}), we have for large $y > 0$,
\begin{equation}
\begin{split}
& \frac{1}{\sqrt{2\pi r_j}}\big\langle R_{free(j)}(k^2 + i0)f_j,e^{ in x/r_j}\big\rangle_j \\
&=  
\left\{
\begin{split}
&y^{1/2}K_{-ik}(|n|y/d_j)\int_0^{\infty}(y')^{1/2}I_{-ik}(|n|y'/d_j)f_{j,n}(y')\frac{dy'}{(y')^2}, \quad n \neq 0, \\
& \frac{1}{2ik}y^{1/2+ik}\int_0^{\infty}(y')^{1/2-ik}f_{j,0}(y')\frac{dy'}{(y')^2}, \quad n = 0.
\end{split}
\right.
\end{split}
\end{equation}
Using $K_{-ik}(z) = K_{ik}(z)$, and noting that
$$
u - u_j^{(in)} = - \sum_n\big\langle R_{free(j)}(k^2+i0)f_j,\frac{e^{inx/r_j}}{\sqrt{2\pi r_j}}\big\rangle_j\,\frac{e^{inx/r_j}}{\sqrt{2\pi r_j}},
$$
we   prove (1). \qed

\medskip
Given $u_j^{(in)}$, $j=1,\cdots,\mu$, one can compute $b_{j,n}$ by observing the asymptotic behavior of $u - u^{(in)}$ in a neighborhood of the cusp. With this in mind, we make the following definition.


\begin{definition}
We call the operator
$$
{\bf S}(k) : {\bf A}_{\infty} \ni \psi^{(in)} \to \psi^{(out)} \in {\bf A}_{-\infty}
$$
the {\it generalized S-matrix}.
\end{definition}


\subsection{Splitting the manifold}
We take a compact submanifold, $\Gamma \subset \mathcal M$, of codimension 1, and split 
$\mathcal M$  into 2 parts, $\mathcal M_{ext}$ and $\mathcal M_{int}$, in the following way:
$$
\mathcal M = \mathcal M_{ext}\cup\mathcal M_{int}, \quad \mathcal M_{ext}\cap\mathcal M_{int} = \Gamma.
$$
Here $\mathcal M_{ext} \setminus \Gamma$ and $\mathcal M_{int} \setminus \Gamma$ 
are assumed to be open submanifolds of 
$\mathcal M$ with boundary $\Gamma$ inheriting the Riemannian structure of $\mathcal M$. 
Assume 
also that $\mathcal M_{ext}$ is non-compact, has infinity common to $\mathcal M_1$ 
and 
 no other infinity.
Recall that the end $\mathcal M_1$ has a cusp. We also assume that $\mathcal M_{sing}$ 
is in the interior of $\mathcal M_{int}$.

Let $-\Delta_g$ be the Laplace-Beltrami operator on $\mathcal M$, $H_{ext}$ and $H_{int}$ 
be $- \Delta_g - 1/4$ defined on $\mathcal M_{ext}$, $\mathcal M_{int}$ with Neumann 
boundary condition on $\Gamma$, respectively. If $\mathcal M$ has only one end 
(i.e. N = 1), $\mathcal M_{int}$ is a compact manifold, and $H_{int}$ has a discrete 
spectrum. If $N \geq 2$, both of $\mathcal M_{int}$ and $\mathcal M_{ext}$ are non-compact, 
and,
{\ntekst although now $\p \mathcal M_{ext}= \p \mathcal M_{int}=\Gamma \neq
\emptyset$, }  the theorems in \S 3 and \S 4 also hold for 
$H_{ext}$, $H_{int}$.
 We denote the inner product of $L^2(\Gamma)$ by
$$
\langle f,g\rangle_{\Gamma} = \int_{\Gamma}f\overline{g}\,dl.
$$

We put  
\begin{equation}
\phi_{n,free} = \left\{
\begin{split}
& y^{1/2-ik}, \quad n = 0,\\
& e^{ inx/r_1}y^{1/2}I_{-ik}(|n|y/r_1), \quad n \neq 0,
\end{split}
\right.
\nonumber
\end{equation}
\begin{equation}
g_n = (H-k^2)\chi_1\phi_{n,free} = [H_{free(1)},\chi_1]\phi_{n,free},
\nonumber
\end{equation}
\begin{equation}
\phi_n^{(+)} =  \chi_1\phi_{n,free} - R(k^2 + i0)g_n.
\label{S4phin}
\end{equation}


\begin{lemma}
We take $\Gamma =\Gamma_0= \{y = y_0\} \subset \mathcal M_1$, $y_0>2$.
Let $k > 0$ and  $k^2 \not\in \sigma_p(H)\cap\sigma_p(H_{int})$. Let $f \in L^2(\Gamma_0)$ satisfy
\begin{equation}
\langle f, \partial_{\nu}\phi_{n}^{(+)}\rangle_{\Gamma_0} = 0  \quad
\forall n \in {\bf Z},
\label{S4forthdelyphi}
\end{equation}
where $\nu$ is the unit normal to $\Gamma_0,\, \p_\nu=\p_y$.
Then $f = 0$.
\end{lemma}

\noindent
{\bf Proof.}  Note that $\Gamma_0$ is naturally identified with $S^{r_1}$.
We define an operator $\delta_{\Gamma_0}' \in {\bf B}(H^{-1/2}(\Gamma_0);H^{-2}(\mathcal M))$ by
\beq \label{9.3.1}
(\delta_{\Gamma_0}'v,w) = \langle v,\partial_{\nu}w\rangle_{\Gamma_0}, \quad 
\forall v \in H^{-1/2}(\Gamma_0), \quad 
\forall w \in H^2(\mathcal M),
\eeq
and define $u = R(k^2 - i0)\delta_{\Gamma_0}'f$ by duality, i.e. for $w \in L^{2,s}$, $s > 1/2$,
\begin{eqnarray*}
(R(k^2 - i0)\delta_{\Gamma_0}'f,w) &=& (\delta_{\Gamma_0}'f,R(k^2 + i0)w) \\
&=& \langle f,\partial_yR(k^2 + i0)w\rangle_{\Gamma_0}. 
\end{eqnarray*}
Note that  
$(H - k^2)u = \delta_{\Gamma_0}'f$ in the sense of distribution, hence, in the classical sense,
\begin{equation}
(H- k^2)u = 0, \quad \mathcal M\setminus\Gamma_0.
\label{S4zerooutsidegammma}
\end{equation}
Considering $H_{free(1)}$ on $M_1 \times (1, \infty)$, we have
\begin{equation}
R_{free(1)}(k^2 -i0)\delta_{\Gamma_0}'f = \frac{1}{2\pi r_1}
\sum_{n\in{\bf Z}}A_n(y)\widehat f_n e^{inx/r_1},\quad f_n=\langle f, e^{inx/r_1} \rangle_1.
\label{S4R0Ndeltagamma}
\end{equation}
Here, taking account of (\ref{S3Rfreek2+i0}), for $n \neq 0$,
\begin{equation}
A_n(y) = \left\{
\begin{split}
\left(y^{1/2}K_{ik}(|n|y/r_1)\right)'\Big|_{y=y_0}y^{1/2}I_{ik}(|n|y/r_1), 
\quad 
y < y_0, \\
\left(y^{1/2}I_{ik}(|n|y/r_1)\right)'\Big|_{y=y_0}y^{1/2}K_{ik}(|n|y/r_1), 
\quad 
y > y_0, 
\end{split}
\right.
\label{S4An}
\end{equation}
and, for $n = 0$,
\begin{equation}
A_0(y) = \left\{
\begin{split}
\left(y^{1/2-ik}\right)'y_0^{1/2+ik}, 
\quad 
y < y_0, \\
\left(y^{1/2+ik}\right)'y_0^{1/2-ik}, 
\quad 
y > y_0.
\end{split}
\right.
\label{S4A0}
\end{equation}
By  (\ref{S4chijresolvent}), (\ref{S4phin}) and (\ref{9.3.1}), when $y>y_0$  and $n \neq 0$, $u = R(k^2 - i0)\delta_{\Gamma_0}'f$ satisfies
\beq \nonumber
& &2\pi r_1\langle u(\cdot, y), e^{inx/r_1}\rangle_{1} 
\\ \nonumber
& &=
y^{1/2}K_{ik}(|n|y/r_1)\int_0^{\infty}\int_{0}^{2\pi r_1}e^{- inx/r_1}(y')^{1/2}I_{ik}(|n|y'/r_1) 
\\ \nonumber
& & \hskip 2cm \times 
 \left\{\chi_1 + [H_{free(1)},\chi_1]R(k^2 - i0)\right\}\delta_{\Gamma}'f\frac{dxdy'}{(y')^2} 
\\ \nonumber
& &= y^{1/2}K_{ik}(|n|y/r_1)\left(\delta_{\Gamma}'f,\left\{\chi_1 - R(k^2 + i0)[H_{free(1)},\chi_1]\right\}\phi_{n,free}\right) 
\\ \nonumber
& &= y^{1/2}K_{ik}(|n|y/r_1)\langle f,\partial_y\phi_n^{(+)}\rangle_{\Gamma_0} = 0.
\eeq
Similarly, one can show that, for large $y$,
$$
\langle u(\cdot, y), 1\rangle_1 = 0.
$$
Therefore, $u = 0$ when $y$ is large enough. Since $(H - k^2)u = 0$ in $\mathcal M_{ext}$, the unique continuation theorem imply that $u = 0$ in $\mathcal M_{ext}$. 
Let $\xi(y) \in C^\infty_0(1, \infty)$ have value $1$ in a neighborhood of $y=y_0$. Then,
$$
R(k^2 - i0)\delta_{\Gamma_0}'f-\xi(y)R_{free(1)}(k^2 - i0)\delta_{\Gamma_0}'f \in C^\infty({\mathcal M}).
$$
Thus, using formulas 
(\ref{S4R0Ndeltagamma}) $\sim$ (\ref{S4A0}), we see that $\partial_yR(k^2 - i0)\delta_{\Gamma_0}'f$ is continuous across $\Gamma_0$. Therefore, in $\mathcal M_{int}$, $u$ satisfies $(H_{int} - k^2)u = 0$ and the Neumann boundary condition on $\Gamma_0$, hence $u = 0$ in $\mathcal M_{int}$. This follows from the assumption $k^2 \not\in \sigma_p(H_{int})$ when $\mathcal M_{int}$ is  compact, and from Theorem 3.8 when $\mathcal M_{int}$ is non-compact. 
We then have $u = 0$ in $\mathcal M$, which implies $\delta_{\Gamma_0}'f = 0$. 
Thus, by (\ref{9.3.1}), $\langle f,\partial_yw\rangle_{\Gamma_0} = 0, \ \forall w \in H^2(\mathcal M)$, which proves  $f = 0$. \qed

\medskip
The generalized S-matrix ${\bf S}(k)$ is an operator-valued $N\times N$ matrix. Let ${\bf S}_{11}(k)$ be its $(1,1)$ entry. For ${\bf a}  \in l^{2,\infty}$, we put
${\bf b} = {\bf S}_{11}(k){\bf a} \in l^{2,-\infty}$, and
$$
\Phi = \sum_{n\in{\bf Z}}a_n\phi_n^{(+)}.
$$
Then, $(H-k^2)\Phi = 0$ and by (\ref{S4uj-}) and (\ref{S4uj+}) it takes the form
$$
\Phi = u_1^{(in)} - u_1^{(out)}
$$
In particular, in $\mathcal M_1$,
$$
u_1^{(in)} = a_{0}y^{1/2-ik} + \sum_{n\neq 0}a_ne^{ inx/r_1}y^{1/2}I_{-ik}(|n|y/r_1),
$$
$$
u_1^{(out)} = b_{0}y^{1/2+ik} + \sum_{n\neq 0}b_ne^{ inx/r_1}y^{1/2}K_{ik}(|n|y/r_1).
$$
Therefore, the knowledge of ${\bf S}_{11}(k)$ is equivalent to 
the observation, for any incoming exponentially growing wave $u_1^{(in)}$ at the 
cusp $\mathcal M_1$, the corresponding outgoing  exponentially decaying wave $u_1^{(out)}$ 
at $\mathcal M_1$.


\subsection{Gel'fand problem, BSP and N-D map}
Before going to proceed, let us recall the {\it Gel'fand inverse boundary-spectral problem}. 
Let $\Omega$ be a compact 
Riemannian manifold with boundary $\p \Omega, \,\Gamma \subset \p \Omega$  
be an open subset, and $- \Delta_g$ be the 
associated Laplace-Beltrami operator. Let $0=\lambda_1 < \lambda_2 < \cdots$ be its Neumann 
eigenvalues without counting multiplicities, and $\varphi_{n,1},\cdots,\varphi_{n,m(n)}$ 
be the orthonormal system of eigenvectors associated with the eigenvalue $\lambda_n$. Let us 
call the set
$$
\Big\{(\lambda_n,\varphi_{n,1}\big|_{\Gamma},\cdots,\varphi_{n,m(n)}\big|_{\Gamma})\Big\}_{n=1}^{\infty}
$$
the boundary spectral data ({\bf BSD}).
The problem raised by Gel'fand  is :

\noindent {\it Do BSD determine the Riemannian manifold $\Omega$?}
This problem was solved by Belishev-Kurylev \cite{BeKu92} using the boundary control method 
(BC-method)   first proposed by
Belishev \cite{Be87} for inverse problems in Euclidean domains.
Later, the method has been developed to study
inverse problems on compact Riemannian manifolds,
\cite{AKKLT, Ku94, KaKu98,KL00,KKL01,KL02,KKL08,LaOk10} and non-compact manifolds
\cite{BKLS,IKL09}. 
The BC-method combines the control theory obtained
from unique continuation results \cite{Ta1,Ta3} with Blagovestchenskii's identity 
that
gives the inner product of the solutions of the wave equation in terms of the 
boundary data.
This identity was originally used in the study of one-dimensional inverse problems,
see \cite{Bl1,Bl3}.

Although it is 
formulated in terms of BSD, what is actually used in the BC-method is the boundary spectral 
projection ({\bf BSP}) defined by
\begin{equation}
\Big\{(\lambda_n, \sum_{j=1}^{m(n)}\varphi_{n,j}(x)\overline{\varphi_{n,j}(y)}\big|_{(x,y)\in\Gamma\times\Gamma})\Big\}_{n=1}^{\infty}.
\label{S4BSPdiscrete}
\end{equation}
This appears in the kernel of the Neumann to Dirichlet map (N-D map)
\begin{equation}
\Lambda(z) : f \to u,
\label{S4NDmap}
\end{equation}
where $u$ is the solution to the Neumann problem
\begin{equation}
\left\{
\begin{split}
& (- \Delta_g - z)u = 0 \quad {\rm in} \quad \Omega, \\
& \partial_{\nu}u = f \in H^{1/2}(\Gamma),
\end{split}
\right.
\label{S4Neumannproblem}
\end{equation}
$\nu$ being the outer unit normal to $\Gamma$, $z \not\in \sigma(-\Delta_g)$.
The N-D map is related to the resolvent $(- \Delta_g - z)^{-1}$ in the following way :
\begin{equation}
\Lambda(z) = \delta_{\Gamma}^{\ast}(- \Delta_g - z)^{-1}\delta_{\Gamma}, \quad 
z \not\in \sigma(-\Delta_g).
\label{S4NDreslvent}
\end{equation}
Here $\delta_{\Gamma} \in {\bf B}\left((H^{1/2}(\Gamma))';\,(H^{1}(\Omega))'\right)$ is the adjoint of 
the trace operator,
$$
r_{\Gamma} : H^1(\Omega) \ni w \to w\big|_{\Gamma} \in H^{1/2}(\Gamma),
$$
\begin{equation}
(\delta_{\Gamma}f,w)_{L^2(\Omega)} = (f,r_{\Gamma}w)_{L^2(\Gamma)}, 
\quad f \in H^{-1/2}(\Gamma), \quad w \in H^1(\Omega),
\label{S4deltagamma}
\end{equation}
and 
we denote by $\left(H^s\right)'$ the dual to $H^s$ with respect to the $L^2-$pairing.
More precisely, we have

\begin{lemma}
To give BSP is equivalent to give the N-D map $\Lambda(z)$ for all $z \not\in \sigma(-\Delta_g)$.
\end{lemma}
We refer for analogous equivalence results for different kind of boundary data to \cite{KL02,KKL04b}.

Let $\Omega$ be non-compact with  asymptotically hyperbolic ends of the type discussed in this paper.
Let $\Gamma \subset \partial\Omega$ be compact, and consider 
its shifted Laplace-Beltrami operator  with Neumann boundary condition, $H= -\Delta_g-\frac14$. 
It has  
continuous spectrum $\sigma_c( H) = [0,\infty)$, and, furthermore, 
$ H$ has a spectral 
representation ${\mathcal F}$ like the one discussed 
in \S3. 
In this case we define the BSP to be the collection
\begin{equation}
\Big\{\delta_{\Gamma}^{\ast}\mathcal F(k)^{\ast}\mathcal F(k)\delta_{\Gamma}\, ; 
\, k > 0\Big\}
\cup\Big\{(\lambda_n,\delta_{\Gamma}^{\ast}P_{n}\delta_{\Gamma})\Big\}_{n=1}^{m}.
\nonumber
\end{equation}
Here $\lambda_n$ is the eigenvalue of $H$, $P_n$ is the associated eigenprojection 
and $m$ is the number of eigenvalues which, in principle could be infinite, see 
Theorem \ref{th.3.6}.
In this case, we extend the N-D map $\Lambda(z)$ for $z  \in \C \setminus\sigma(H)$ 
by using the 
solution $u$ of (\ref{S4Neumannproblem}). Note that we can extend the
definition of $\Lambda(z)$ for $z =k^2 \pm i0\notin \sigma_p(H)$ 
by using  the outgoing or incoming radiation conditions. 
Then Lemma 4.6 also holds in this case. 
(See \cite{IsKu10}, Chap. 5, \S 3 and \S 4, \cite{IKL09}, Lemma 5.6.)

Denote by $G(z; X, Y),\, z \in \C \setminus\sigma(H)$, the 
Schwartz kernel of $(H-z)^{-1}$.
Since
\ba
(H-z)^{-1}= \sum_{n=1}^m \frac 1{\lambda_n-z}P_n+
\int_{0}^\infty \, \frac 1{k^2-z}\mathcal F(k)^{\ast}\mathcal F(k)\,dk,
\ea
we have, in view of (\ref{S4NDreslvent}), that
\beq
\hspace{1cm}
G(z;\, \cdotp, \cdotp)\Big|_{\Gamma\times\Gamma} =
\sum_{n=1}^m \frac 1{\lambda_n-z}\delta_{\Gamma}^{\ast}P_n\delta_{\Gamma}
+\int_0^\infty\,  \frac 1{k^2-z}\delta_{\Gamma}^{\ast}
\mathcal F(k)^{\ast}\mathcal F(k)\delta_{\Gamma}\, dk.\hspace{-1cm}
\label{decom of f}
\eeq
Here the left-hand side is understood as the Schwartz kernel of the operator in the right-hand side of the formula.


\subsection{Generalized S-matrix and N-D map}
Returning to our problem concerning  2-dimensional non-compact surfaces with conical singularities, we take $\Omega = \mathcal M_{int}$ with 
$\Gamma_0 = \{X \in  \mathcal M_1:\, y=y_0\}$, $y_0 >2$. We define the N-D map for $\mathcal M_{int}$ by (\ref{S4NDmap}) and (\ref{S4Neumannproblem}). 

Now suppose we are given two manifolds $\mathcal M^{(i)}$, $i = 1,2$, satisfying the assumptions 
(A-1) $\sim$ (A-4) in \S 1. Let $-\Delta^{(i)}$ be the Laplace-Betrami operator of 
$\mathcal M^{(i)}$. Assume that $\mathcal M^{(i)}$ has $N_i$ numbers of ends, and 
let ${\bf S}^{(i)}_{11}(k)$ be the $(1,1)$ entry of the generalized S-matrix for 
$H^{(i)}=-\Delta^{(i)}-\frac14$. 

Assuming that $r_1^{(1)}=r_1^{(2)}$,  we can naturally identify $\mathcal M_{1}^{(1)}$ 
and $\mathcal M_{1}^{(2)}$. Taking $\Gamma_0$ as above, we split 
$\mathcal M^{(i)}$ into $\mathcal M_{int}^{(i)} \cup \mathcal M_{ext}^{(i)}$ 
by using 
$\Gamma_0$. 
Let $H^{(i)}_{int}=-\Delta_{int}^{(i)}-\frac14$ be the shifted 
Laplace-Beltrami operator of $\mathcal M_{int}^{(i)}$ 
with Neumann boundary condition on $\Gamma_0$, and define the N-D map $\Lambda^{(i)}(z)$ 
for $H_{int}^{(i)}$. With this preparation, we can prove the following lemma.


\begin{lemma} Let $k > 0$, $k^2 \not\in \sigma_p(H^{(1)}) \cup\sigma_p(H^{(2)})\cup\sigma_p(H_{int}^{(1)})\cup\sigma_p(H^{(2)}_{int})$.
If ${\bf S}^{(1)}_{11}(k) = {\bf S}^{(2)}_{11}(k)$, we have 
$\Lambda^{(1)}(k^2 + i0) = \Lambda^{(2)}(k^2 +i0).$
\end{lemma}

\noindent
{\bf Proof.}
For $i = 1, 2$, we construct $\phi_n^{(i)}$ as in (\ref{S4phin}), and put $u = \phi_n^{(1)} - \phi_n^{(2)}$. Then $u$ satisfies 
\ba
(H^{(i)}-k^2)u=(-\Delta^{(i)} - k^2 -\frac14)u = 0 \quad \hbox{for} \,\,
X \in\mathcal M_{ext}^{(1)} = \mathcal M_{ext}^{(2)},\,\, y(X)>2,
\ea
due to $\chi_1 = 1$ there, see (\ref{2.3.1}). Since ${\bf S}^{(1)}_{11}(k) = {\bf S}^{(2)}_{11}(k)$, by the same argument as in the proof of Lemma 4.5, we have $u = 0$ for 
$X \in\mathcal M_{ext}^{(1)} = \mathcal M_{ext}^{(2)},\,\, y(X)>2$. 
Hence, $\partial_{\nu}\phi_n^{(1)} = \partial_{\nu}\phi_n^{(2)}$ on $\Gamma_0$. 

In $\mathcal M_{int}^{(i)}$, $\phi_n^{(i)}$ is the outgoing solution of the equation 
$(H^{(i)}_{int} - k^2) \phi_n^{(i)} = 0$. Hence,  
$\partial_{\nu}\phi_n^{(i)}\big|_{\Gamma_0} = \Lambda^{(i)}(k^2+i0)\phi_n^{(i)}\big|_{\Gamma_0}$,
where we again use that $\chi_1 =1$ near $\Gamma_0$. This implies
\begin{equation}
\Lambda^{(1)}(k^2+i0)\phi_n^{(1)}\big|_{\Gamma_0} = 
\Lambda^{(2)}(k^2+i0)\phi_n^{(2)}\big|_{\Gamma_0}, \quad \forall n.
\label{S4NDcoincide}
\end{equation} 
Lemma 4.5 implies that the linear span of $\{\partial_{\nu}\mathcal \phi_n^{(i)}\big|_{\Gamma_0}\, ; \, n \in {\bf Z}\}$ is dense in $L^2({\Gamma_0})$. Therefore, 
by (\ref{S4NDcoincide}), $\Lambda^{(1)}(k^2+i0) = 
\Lambda^{(2)}(k^2+i0)$. \qed


\begin{cor} \label{Cor. 4.8} Let $(a,b)$ be an interval such that $(a ,b) \cap 
\left(\sigma_p(H^{(1)})\cup\sigma_p(H^{(2)})  \right)= \emptyset$, and assume that ${\bf S}^{(1)}_{11}(k) = {\bf S}^{(2)}_{11}(k)$ for $k^2 \in (a, b)$.
 Then $\Lambda^{(1)}(z)=\Lambda^{(2)}(z)$ if 
$ z \notin \sigma(H_{int}^{(1)})\cup\sigma(H_{int}^{(2)})$.
 Moreover, BSP's for $H_{int}^{(1)}$ and $H_{int}^{(2)}$ and Green's kernels $G^{(i)}(z; X, Y)$ for $(H^{(i)} -z)^{-1}$, $i = 1, 2$, coincide on $\Gamma_0\times\Gamma_0$.
\end{cor}
\noindent
{\bf Proof.} For $f \in H^{1/2}(\Gamma_0),$ let $ F \in H^2_0({\mathcal M}^{(i)}_{int})$ satisfy
\begin{equation} \label{E1}
\p_{\nu} F|_{\Gamma_0}= f, \quad \hbox{supp}(F) \subset S^{r_1} \times [2, y_0].
\end{equation}
Then,
$$
\Lambda^{(i)}(z) f = r_{\Gamma_0}\left( F - (-\Delta^{(i)} -z)^{-1} w_f(z) \right),
$$
where $w_f(z)= -(H_{int}^{(i)} -z) F$ is independent of $i=1, 2,$ due to (\ref{E1}).
Note that $\Lambda^{(i)}(z) f$ is analytic, if $z \notin \sigma(H^{(i)}_{int})$, and have a limit,
$\Lambda^{(i)}(k^2 \pm i 0)f$,
when $z \to k^2 \pm i0,\, k^2 \notin \sigma_p(H^{(i)}_{int})$.
Using Lemmas 4.7, 4.6 and (\ref{decom of f}) we obtain the result. \qed


\section{Uniqueness of inverse scattering}


\subsection{Blagovestchenskii's identity}

To prove the uniqueness of the inverse scattering problem we start with
some auxiliary results.
Let $\Omega$ be a (possibly non-compact) Riemannian surface with conical singularities 
(and asymptotically hyperbolic ends)
and $H=-\Delta_g-\frac 14$ be the Hamiltonian
corresponding to Neumann boundary condition on $\p \Omega$. 
We denote by ${\mathcal F}(k)$ the Fourier transform associated to $H$ and by
$P_j$ the orthogonal projections corresponding to eigenvalues $\lambda_j$ of $H$
using the convention that  ${\mathcal F}(k)=0$ when $\Omega$ is compact.
Let $\Gamma\subset \p\Omega$ be open. 
Consider the solution $u^f(X,t)$ of the initial boundary value problem
\begin{equation}
\left\{
\begin{split}
& \partial_t^2u - \Delta_g u-\frac14 u=0, \quad {\rm in} \quad \Omega \times \R_+, \\
& u\big|_{t=0} = \partial_tu\big|_{t=0} = 0, \quad {\rm in } \quad \Omega,\\
& \partial_{\nu}u = f, \quad {\rm in} \quad \partial\Omega \times \R_+, \quad {\rm supp}\,f \subset \Gamma \times{\R}_+.
\end{split}
\right.
\label{S5:IBVPWave}
\end{equation}

Let
$$
B(t,\lambda) = \left\{\begin{array}{cl}\frac{\sin(\sqrt{\lambda}t)}{\sqrt\lambda}& \hbox{for }\lambda\not =0,\\
t&\hbox{for } \lambda =0.\end{array}\right.
$$


\begin{lemma} \label{lem: step 2}
Assume that we are given the curve  $\Gamma \subset \p \Omega$, the length element $dl$ on $\Gamma$ and the BSP of $H$ on $\Gamma$. 
Then, for any given $f,h\in C^\infty_0(\Gamma\times \R_+)$ and $t,s>0$, these data uniquely determine
\ba
(u^f(t),u^h(s))_{L^2(\Omega)}=\int_{\Omega} u^f(X,t)\,\overline {u^h(X,s)}\,dS_X
\ea
and
\ba
(u^f(t),1)_{L^2(\Omega)}=\int_{\Omega} u^f(X,t)\,dS_X.
\ea
Moreover, the hyperbolic N-D map $R^T_\Gamma:f\mapsto u^f|_{ \Gamma \times (0,T)}$
can be written in terms of BSP as
\beq\nonumber \quad \quad
R^T_\Gamma f(\,\cdotp,t)&=&\int_0^t dt'\bigg( 
\sum_{n=1}^m B(t-t',\lambda_n)\delta_{\Gamma}^{\ast}P_n\delta_{\Gamma} +
\\  \label{decom of R map} & &\quad+
\int_0^\infty dk\, B(t-t',k^2)\delta_{\Gamma}^{\ast}\mathcal F(k)^{\ast}\mathcal F(k)\delta_{\Gamma}\bigg)f(\,\cdotp,t').
\eeq 
\end{lemma} 

\noindent {\bf Proof.} 
The solution $u^f(t)$ can be  written as 
\ba
& &u^f(X,t) = 
\\ \nonumber
& &\int_0^t dt'\bigg( 
\sum_{n=1}^m B(t-t',\lambda_n)P_n\delta_{\Gamma} 
+
\int_0^\infty dk\, B(t-t',k^2)\mathcal F(k)^{\ast}
\mathcal F(k)\delta_{\Gamma}\bigg)f(\,\cdotp,t').
\ea
Restricting this equation to $\Gamma$, we  prove (\ref{decom of R map}).

Using the similar decomposition for $u^h(s)$,  we  obtain the following formula:
\beq \label{eq:Sect5Blagovecont}
& &\hspace{1cm}(u^f(t),u^h(s))_{L^2}^2= \\
&&\hspace{-1cm} \int_0^tdt'\int_0^sds' \int_{\Gamma}  dl_X  \int_{\Gamma}  dl_Y \,\widetilde K(t-t',s-s',X,Y)
f(X,t')h(Y,s').\hspace{-1cm}
\nonumber
\eeq
Here 
\begin{equation}
\begin{split}
\widetilde K(t,s, \cdot, \cdot)= &
\sum_{n=1}^m B(t,\lambda_n) B(s, \la_n) \delta_{\Gamma}^{\ast}P_n\delta_{\Gamma} \\
& + \int_0^\infty dk\, B(t,k^2) B(s, k^2) \delta_{\Gamma}^{\ast}\mathcal F(k)^{\ast}\mathcal F(k)\delta_{\Gamma},
\end{split}
\nonumber
\end{equation}
where the left-hand side is understood as the Schwartz kernel of the operator in the 
right-hand side.
 
Moreover, as   $J(t)=(u^f(t),1)_{L^2}$ satisfies the differential equation 
\ba
\p_t^2 J(t)=(\p_t^2 u^f(t),1)_{L^2}
=(\Delta_g u^f(t),1)_{L^2}
=\int_{\Gamma}  f(Y,t) dl_Y
\ea
and initial conditions $J(0)=\p_t J(t)|_{t=0}=0$, we see that
  \ba
& &(u^f(t),1)_{L^2}
=
\int_\Gamma dl_Y \int_0^t dt' \,B(t-t',0)f(Y,t').
\ea
\qed

\medskip
Above, the formula (\ref{eq:Sect5Blagovecont}) is a generalization of Blagovestchenskii identity (see 
\cite[Theorem 3.7]{KKL01}) for Riemannian surfaces with conic singularities. 




Next we will apply these formulas to compute the area of the domain of influence
\begin{equation}
\Omega(\tilde \Gamma,T)=\{X \in \Omega:\, d_g(X, \tilde \Gamma) \leq T   \},
\quad \tilde \Gamma \subset  \Gamma,
\label{S5Domainofinfluence}
\end{equation}
where $d_g$ denotes the distance in $\Omega$ with respect to $g$. We denote the area of $\Omega(\tilde \Gamma,T)$ by $S_g(\Omega(\tilde \Gamma,T))$.


\begin{lemma} \label{lem: volumes}
Assume that we are given the curve  $\Gamma$, the length element $dl$  on $\Gamma$ and the 
BSP of $H$ on $\Gamma$. 
Then, for any given open set $\tilde \Gamma \subset \Gamma$ and $T>0$, these data 
uniquely 
determine $S_g(\Omega(\tilde \Gamma,T))$.
\end{lemma}

\noindent {\bf Proof.} Let $w\in L^2(\Omega)$ be a function such that $w=1$ in 
$\Omega(\tilde \Gamma,T)$.
For $f\in  C^\infty_0(\tilde \Gamma\times (0,T))$, real-valued, we define the quadratic functional 
\ba
I_T(f)&=& \|u^f(\cdot, T)-w\|^2_{L^2(\Omega)}-\|w\|^2_{L^2(\Omega)}.
\ea
Since
 $
 \supp(u^f(\cdot, T))\subset \Omega(\tilde \Gamma,T), 
$
 we have
\beq \label{5.5a}
I_T(f)
&=& \|u^f(\cdot, T)\|^2_{L^2(\Omega)}-2 (u^f(\cdot, T),1)_{L^2(\Omega)}.
\eeq
Hence, by Lemma \ref{lem: step 2},   we can compute $I_T(f)$ for any 
$f\in C^\infty_0(\tilde \Gamma\times (0,T))$ uniquely by using BSP and $dl$ on $\Gamma$. In the sequel, this is phrased as {\it we can compute}.

Now we use again the fact that, for $f\in  C^\infty_0(\tilde \Gamma\times (0,T))$,
 $\supp(u^f(\cdot, T))\subset \Omega(\tilde \Gamma,T)$ so that (\ref{5.5a}) yields that
\ba
I_T(f)&=& \|u^f(\cdot, T)-\chi_{(\tilde \Gamma,T)}\|^2_{L^2( \Omega)}- 
\|\chi_{(\tilde \Gamma,T)}\|^2_{L^2( \Omega)},
\ea
where $\chi_{(\tilde \Gamma,T)}$ is the characteristic function of $\Omega(\tilde \Gamma,T)$.
Thus, 
\beq\label{minimimum}
I_T(f)\geq - S_g(\Omega(\tilde \Gamma,T)), \quad\hbox{for all }f\in  C^\infty_0(\tilde \Gamma\times (0,T)).
\eeq
By Tataru's controllability theorem, see \cite{Ta1} and e.g.\ \cite{KKL01},  there is a sequence $h_j\in  C^\infty_0(\tilde \Gamma\times (0,T))$, such that
\ba
\lim_{j\to \infty} u^{h_j}(\cdot, T)=\chi_{\Omega(\tilde \Gamma,T)}\quad \hbox{in }L^2(\Omega).
\ea
For this sequence,
\beq\label{minimizing 0}
\lim_{j \to \infty}I_T(h_j) = - S_g(\Omega(\tilde \Gamma,T)).
\eeq
On the other hand, if 
$f_j \in C^\infty_0(\tilde \Gamma\times (0,T))$ is a minimizing sequence for $I_T$, i.e.,
\beq\label{minimizing I}
\lim_{j\to \infty}I_T(f_j)=m_0:=\inf \{I_T(f);\,  f\in C^\infty_0(\tilde \Gamma\times (0,T))\},
\eeq
then, by  (\ref{minimimum}) and (\ref{minimizing 0}),
\ba
\lim_{j\to \infty} u^{f_j}(\cdot, T)=\chi_{\Omega(\tilde \Gamma,T)}\quad \hbox{in }L^2(\Omega).
\ea
Thus, using any sequence $(f_j)$ satisfying (\ref{minimizing I}), we can compute
\ba
S_g(\Omega(\tilde \Gamma,T))=\lim_{j\to \infty} (u^{f_j}(\cdot, T),u^{f_j}(\cdot, T))_{L^2(\Omega)}.
\ea
\qed

\subsection{Reconstruction near $\Gamma_0$} \label{Green}


To prove Theorem \ref{thm: main result} our first aim is 
to show that $ {\mathcal M}_{reg}^{(1)}$ and $ {\mathcal M}_{reg}^{(2)}$
are isometric. The proof is based on the procedure 
of the continuation of 
Green's functions, $G^{(i)}(z; X, Y),\, i=1,2$,
of the operators $H^{(i)}_{int}$. 

We are going to prove the uniqueness for the inverse problem
 step by step by  constructing relatively open subsets 
 ${\mathcal M}^{(1),rec} \subset {\mathcal M}_{int}^{(1)}$ 
 and  ${\mathcal M}^{(2),rec} \subset {\mathcal M}_{int}^{(2)}$, which are isometric
 and enlarge these sets at each step.
 In the following, when  
 ${\mathcal M}^{(1),rec} \subset {\mathcal M}_{reg}^{(1)}\cap {\mathcal M}_{int}^{(1)}$ 
 and  ${\mathcal M}^{(2),rec} \subset {\mathcal M}_{reg}^{(2)}\cap {\mathcal M}_{int}^{(2)}$  
 are relatively open connected sets and 
 \ba
\Phi^{rec}:\, {\mathcal M}^{(1),rec} \rightarrow {\mathcal M}^{(2),rec},
\ea
is a diffeomorphism, we say that the triple 
  $({\mathcal M}^{(1),rec},{\mathcal M}^{(2),rec},\Phi^{rec})$ is {\it admissible}
  if $\Phi^{rec}:\, {\mathcal M}^{(1),rec} \to {\mathcal M}^{(2),rec}$
 is an isometry, that is, $(\Phi^{rec})_*g^{(1)}=g^{(2)} $
and the values of Green's functions  $G^{(i)}(z, X, Y)$ on
${\mathcal M}^{(i),rec} $ satisfy, for $X, Y \in {\mathcal M}^{(1),rec},$
the relation
\beq\label{eq: Greens functions related}
G^{(2)}(z; \Phi^{rec}(X), \Phi^{rec}(Y))=G^{(1)}(z; X, Y),\quad\hbox{for }
\ z\in \C\setminus \R.
\eeq

First we consider  Green's functions in  the set  
$$
N=\Gamma_0 \times (2, y_0] \subset {\mathcal M}_1.
$$


\begin{lemma}\label{lem. N1 admissible} When $N$ is considered both as  a subset
${\mathcal M}_{int}^{(1)}$ and ${\mathcal M}_{int}^{(2)}$ and $I:N\to N$
is the identity map, then the triple $(N,N,I)$ is admissible.
\end{lemma}

\noindent
{\bf Proof.} 
By the assumpton (A-3) and $r_1^{(1)}=r_1^{(2)}$, the map $I:(N,g^{(1)})\to (N,g^{(2)})$ is an isometry.
By Corollary \ref{Cor. 4.8}, we know that
\ba
G^{(1)}(z; X, Y)=G^{(2)}(z; X, Y),\quad z \in \C \setminus \R_+,\,\, X, Y \in \Gamma_0.
\ea
Let  $z\in \C\setminus \R$.
Since Green's
function $G^{(i)}(z; X, Y)$ satisfies the elliptic equation,
\beq\label{extra 1}
(-\Delta^{(i)} -\frac14 -z)G^{(i)}(z; \cdot, Y)= \delta_Y,\quad\hbox{on } {\mathcal M}^{(i)}_{int},
\\ \nonumber
\p_\nu G^{(i)}(z; \cdot, Y)|_{\Gamma_0}=0, \nonumber
\eeq
and $g^{(1)}(X)=g^{(2)}(X),$ for $ X \in N$,
 we can use the principle 
of unique continuation with respect to $X$  to show that 
$G^{(1)}(z; X, Y)=G^{(2)}(z; X, Y)$ if $X \in N$,
$Y\in \Gamma_0$.
Moreover, as $ G^{(i)}(z; X, Y)= \overline {G^{(i)}(\overline z; Y, X)}$, Green's
function satisfies an elliptic equation analogous to  (\ref{extra 1}) also
in the $Y$ variable. Thus, using  the principle 
of unique continuation with respect to $Y$,  we see that 
$G^{(1)}(z; X, Y)=G^{(2)}(z; X, Y)$ for $X, Y \in N $ and $z\in \C\setminus \R$. \hfill\qed \medskip


\subsection{Continuation by Green's functions} \label{Green2}

To reconstruct subsets of manifolds   ${\mathcal M}^{(i)}_{int}$, $i=1,2$,
by continuing Green's function, we need the the following result telling
 that the values  of Green's functions
identify the points of the manifold. 

	
\begin{lemma} \label{identification} Let  $X_1,X_2\in {\mathcal M}^{(i)}_{int}$ be such 
that
\beq \label{3.3.4}
G^{(i)}(z, X_1, Y)=G^{(i)}(z, X_2, Y)\,  
\eeq
for all $\, Y \in \Gamma_0$ and some $z \in \C \setminus \R$.
Then $X_1=X_2.$
\end{lemma}
\noindent
{\bf Proof.} 
Using the unique continuation principle for the solutions
of elliptic equations as above after (\ref{extra 1}), we see that (\ref{3.3.4}) implies that 
$G^{(i)}(z, X_1, Y)=G^{(i)}(z, X_2, Y)$, for all $Y \in {\mathcal M}_{int}^{(i)}\setminus
\{X_1,X_2\}$.
As the map $Y\mapsto G^{(i)}(z, X, Y)$ is bounded in the compact
subsets of ${\mathcal M}_{int}^{(i)}\setminus\{X\}$ and tends
to infinity as $Y$ approaches $X$, this proves that $X_1=X_2$. 	
\qed

\begin{remark} \label{Remark:3}
Lemma \ref{identification} has the following important
consequence: If the triples $(N_1^{(1)},N_1^{(2)},\Phi_1)$
and $(N_2^{(1)},N_2^{(2)},\Phi_2)$ are
admissible and $N_1^{(1)}\cap N_2^{(1)}\not = \emptyset,
$ then, by  Lemma \ref{identification}, the maps $\Phi_1(x)$ and 
$\Phi_2(x)$ have to coincide  in $N_1^{(1)}\cap N_2^{(1)}$.
Moreover, if $N_3^{(i)}=N_1^{(i)}\cup N_2^{(i)}$, $i=1,2,$ 
and
\beq\label{eq: glue}
\Phi_3(x)=\left\{\begin{array}{cl}
\Phi_1(x), & \hbox{for } x\in N_1^{(1)}, \\
\Phi_2(x), & \hbox{for } x\in N_2^{(1)},
\end{array}
\right.
\eeq
then, by  Lemma \ref{identification}, the map 
$\Phi_3:N_3^{(1)}\to N_3^{(2)}$ is bijective
and hence a diffeomorphims. This implies that the 
triple $(N_3^{(1)},N_3^{(2)},\Phi_3)$ is
admissible
\end{remark}


The procedure of constructing the isometry between $\mathcal M^{(1)}_{int}$ and $\mathcal M^{(2)}_{int}$ consists of extending the admissible triple $(\mathcal M^{(1),rec},\mathcal M^{(2),rec},\Phi^{rec})$. 
In the first step, we apply Lemma \ref{lem. N1 admissible} 
to the triple $(N,N,I)$. In the subsequent steps 
we always assume that $N\subset {\mathcal M}^{(i),rec}$.

{\ntekst Let 
$q_i\in  {\mathcal M}^{(i),rec}$, $i=1,2$, 
\beq \label{eq: basic assumptions}
\hspace{1cm}\Phi^{rec}(q_1)=q_2,\quad d^{(i)}(q_i, \, \Gamma_0) > (y_0-2)/2,
\eeq
where $d^{(i)}$ denotes the distance on $ {\mathcal M}^{(i)}$.
Let $R=R(q)>0$ be sufficiently small so that $R <(y_0-2)/4$ and the Riemannian 
normal coordinates, centered at $q_i$, are well defined in $B^{(i)}(q_i, 2R)$,
i.e. the ball of the radius $2R$ with respect to the distance $d^{(i)}$.
Assume also that $R$ is so small that ${\mathcal O}^{(i)}=B^{(i)}(q_i,R) $ 
satisfy
\beq \label{eq: basic assumptions2}
& &\overline {\mathcal O}^{(i)} \subset {\mathcal M}^{(i),rec}\setminus \Gamma_0.
\eeq
Then  $\Phi^{rec}({\mathcal O}^{(1)})={\mathcal O}^{(2)}$,
${\mathcal M}^{(i),rec}\setminus \overline {\mathcal O}^{(i)}$ are connected
and  ${\mathcal O}^{(i)}$ has smooth boundary.  }


%

%
Denote $\Omega_{\mathcal O}^{(i)}={\mathcal M}_{int}^{(i)}\setminus {\mathcal O}^{(i)}$.
 We put $H_{\mathcal O}^{(i)}= - \Delta^{(i)}-\frac 14$ in $\Omega_{\mathcal O}^{(i)}$ endowed 
 with the Neumann boundary condition:
\begin{equation}
\partial_{\nu}v = 0 \quad {\rm on} \quad \partial\Omega_{\mathcal O}^{(i)},
\label{eq:OmegainBC}
\end{equation}
$\nu$ being the unit normal to the boundary.

Let
$z \in \C\setminus\R$ and consider the Schwartz kernel 
$G_{\O}^{(i)}(z; X,Y)$ of the operator $(H_{\O}^{(i)} -z)^{-1}$. It satisfies the equation
\beq\label{elliptic eqs}
& &\left(-\Delta^{(i)}-\frac 14 -z \right)G_{\O^{(i)}}(z; \cdotp,Y)=\delta_{Y},\quad Y\in \Omega_{\O}^{(i)},\\
& &\nonumber
\p_\nu G_{\O}^{(i)}(z; \cdotp,Y)|_{\Gamma_0\cup \p {\O}^{(i)}}=0.
\eeq

Let ${\mathcal O}^{(i)}\subset  {\mathcal M}^{(i),rec}$, $i=1,2$ be relatively compact subsets
with smooth boundaries (which later will be chosen to be the balls described earlier). Let
 $\Phi:\p \O^{(1)}\to \p \O^{(2)}$ be  a diffeomorphism. 
Let $(\delta_{\O^{(i)}}^{\ast}\mathcal F^{(i)}(k)^{\ast}\mathcal F^{(i)}(k)\delta_{\O^{(i)}})_{k\in\R_+}$
and $(\lambda_n^{(i)})_{n=1}^{m_i}$ and $(P_n^{(i)})_{n=1}^{m_i}$ 
be the BSP related to operator $H_{\O}^{(i)}$ on $\p \O^{(i)}$, $i=1,2$. 
We say that the 
 BSP related to operators $H_{\O^{(1)}}$ on $\p \O^{(1)}$ 
and $H_{\O^{(2)}}$ on $\p \O^{(2)}$ are $\Phi$-{\it related} 
 if $m_1=m_2$ and, 
for all $h\in C^\infty(\p \O^{(2)})$, $k>0$, and $ j=1,2,\dots,m_1$,
 we have 
\ba
& &\delta_{\O^{(1)}}^{\ast}\mathcal F^{(1)}(k)^{\ast}\mathcal F^{(1)}(k)((\Phi^* h)\delta_{\O^{(1)}})\Phi^{\ast}=\Phi^*\bigg(\delta_{\O^{(2)}}^{\ast}\mathcal F^{(2)}(k)^{\ast}\mathcal F^{(2)}(k)(h\delta_{\O^{(2)}})\bigg),\quad\\
& &\lambda_n^{(1)}=\lambda_n^{(2)},\quad
\delta_{\O^{(1)}}^{\ast}P^{(1)}_{n}((\Phi^* h)\delta_{\O^{(1)}})\Phi^{\ast}=\Phi^*\bigg(\delta_{\O^{(2)}}^{\ast}P^{(2)}_{n}(h\delta_{\O^{(2)}})\bigg).
\ea
Note that $\Phi^{\ast}$ induces a bounded operator $: H^s(\partial\mathcal O^{(2)}) \to H^s(\partial\mathcal O^{(1)})$, which is denoted by $\Phi^{\ast}$ again.

\begin{lemma} \label{lem. determination of Green} 
Let $({\mathcal M}^{(1),rec},{\mathcal M}^{(2),rec},\Phi^{rec})$
be an admissible triple and ${\mathcal O}^{(i)}$, $i=1,2$ be relatively compact subsets
of $ {\mathcal M}^{(i),rec}$ such that  ${\mathcal O}^{(2)}= \Phi^{rec}({\mathcal O}^{(1)})$
and 
${\mathcal M}^{(i),rec} \setminus \overline {\mathcal O}^{(i)}$
are connected.
Let  $ G_{\O}^{(i)}(z; X,Y)$, $z\in \C\setminus \R$ be the 
Schwartz kernels of $ (H_{\O}^{(i)}-z)^{-1}$. Then
\beq\label{eq: two greens functions}
\hspace{1cm}G_{\O}^{(1)}(z; X,Y)=G_{\O}^{(2)}(z;\Phi^{rec}( X),\Phi^{rec}(Y)),\quad X, Y \in {\mathcal M}^{(1),rec}.\hspace{-1cm}
\eeq
Moreover, the BSPs related to operators $H_{\O}^{(1)}$ on $\p \O^{(1)}$ 
and $H_{\O}^{(2)}$ on $\p \O^{(2)}$ are $\Phi^{rec}$-related.
\end{lemma}
\noindent
{\bf Proof.} 
Skipping for a while the superscript $^{(i)}$, we start the proof by 
assuming that
we are given $G(z; X, Y)$ for $X,Y\in  {\mathcal M}^{rec}\subset {\mathcal M}_{int}$
 and $z\in \C\setminus \R$ and showing
that if
${\mathcal O}$ is a relatively compact subset with a smooth boundary
of the open set $ {\mathcal M}^{rec}$ such that
${\mathcal M}^{rec} \setminus \overline {\mathcal O}$ is connected, then
we can determine $ G_{\O}(z; X,Y)$ 
for  $X,Y\in{\mathcal M}^{rec} \setminus \overline {\mathcal O}$ and $z\in \C\setminus \R$.


{\nntext 
To show this, let us denote  by $G^{ext}_{\mathcal O}(z; X,Y)$ some smooth extension of 
$X\mapsto G_{\mathcal O}(z; X,Y)$ into $\mathcal O$, where
 $Y\in {\mathcal M}^{rec} \setminus \overline {\mathcal O}$. Then 
$$
\left(-\Delta-\frac 14 -z \right)G^{ext}_{\mathcal O}(z; \cdotp,Y) - \delta(\cdotp,Y) = 
F(\cdotp,Y) \in C^{\infty}({\mathcal M}^{rec}),
$$
where ${\rm supp}\, F(\cdot,Y) \subset {\overline{\mathcal O}}$ is fixed. 
Therefore,
$$
G_{\mathcal O}(z; X,Y) = G(z;X,Y) + 
\int_{\mathcal O}G(z; X,Y')F(Y',Y)dS_{Y'}.
$$
In particular, due to boundary condition (\ref{elliptic eqs}), if $X \in \partial\mathcal O,$
\begin{equation}
\partial_{\nu(X)}G(z; X,Y) + \int_{\O}\partial_{\nu(X)}G(z; X,Y')F(Y',Y)dS_{Y'} = 0,  
\label{Gznormal}
\end{equation}
where $\nu(X)$ is the unit normal to $\mathcal O$ at $X$. On the other hand, if $F(\cdotp, Y)
 \in C^{\infty}({\mathcal M}^{rec})$, ${\rm supp}\, F(\cdotp, Y) \subset \overline{\mathcal O}$, 
 satisfies the equation (\ref{Gznormal}), then the function
\begin{equation}
G(z; X,Y) + \int_{\O}G(z; X,Y')F(Y',Y)dS_{Y'}, \quad 
X,Y \in {\mathcal M}^{rec}\setminus  \O,
\label{extra}
\end{equation}
is equal to $G_{\mathcal O}(z;X,Y)$. As we have in our disposal $G(z; X,Y)$ for $X, Y \in {\mathcal M}^{rec}$, 
we can verify, for any given function $F$, if it satisfies the equation (\ref{Gznormal}) or not.
As the equation (\ref{Gznormal}) has, for every $Y \in \p \O$, at least one solution, 
this implies that 
we can  find some solution $F$ for the equation (\ref{Gznormal})  and thus determine
the values of $G_\O(z; X, Y)$ for $X, Y\in \p \O$ and $z\in \C\setminus \R$.}

Let $G_\O(z;X,Y)$ be the Schwartz kernel of the operator
  $\delta_{\p \O}^{\ast}(H_\O-z)^{-1}\delta_{\p \O} $. 
 We have 
 \begin{equation}
\begin{split}
\delta_{\p \O}^{\ast}(H_\O-z)^{-1}\delta_{\p \O} = &
\int_0^{\infty}(\lambda-z)^{-1}\delta_{\p \O}^{\ast}{\mathcal F_\O}(\lambda)^{\ast}{\mathcal F_\O}(\lambda)\delta_{\p \O}d\lambda \\
 & \ \ \ \ \ + 
\sum_{n=1}^{m}(\lambda_n-z)^{-1}\delta_{\p \O}^{\ast}P_n\delta_{\p \O}.
\end{split}
\nonumber
\end{equation}
 Using this we see that, for $ \lambda>0$, 
 \beq\label{eq: jump of green's function}
 & &\delta_{\p \O}^{\ast}{\mathcal F_\O}(\lambda)^{\ast}{\mathcal F_\O}(\lambda)\delta_{\p \O}=\\ \nonumber
 & &\frac 1{2\pi i}\lim_{\e\to 0+} 
 \left(\delta_{\p {\mathcal O}}^{\ast}(H_\O-\lambda-i\e)^{-1}\delta_{\p {\mathcal O}} -
 \delta_{\p {\mathcal O}}^{\ast}(H_\O-\lambda+i\e)^{-1}\delta_{\p {\mathcal O}} \right)
 \eeq
 and that $\lambda_n$ are the poles
 of the meromorphic function 
 $\delta_{\p {\mathcal O}}^{\ast}(H_\O-z)^{-1}\delta_{\p {\mathcal O}} $ in 
  $ \C$. Its
  residues  satisfy
  \beq\label{eq: residues}
\hbox{res}_{z=\lambda_n}  \delta_{\p {\mathcal O}}^{\ast}(H_\O-z)^{-1}\delta_{\p {\mathcal O}}
=-\delta_{\p \O}^{\ast}P_n\delta_{\p \O}.
\eeq
 Summarizing the above, we have shown the set ${\mathcal M}^{rec}$ with its metric and  
 values of
 Green's function  $G_\O(z; X, Y)$ for
  $X, Y\in \p\O$ and $z\in \C\setminus \R$  determine the 
BSP on $\p\O$.

As  $({\mathcal M}^{(1),rec},{\mathcal M}^{(2),rec},\Phi^{rec})$ is admissible,
we see that $F^{(2)}$ solves equation (\ref{Gznormal}) on $\mathcal M^{(2),rec} $ 
if and only if  $F^{(1)}=(\Phi^{rec})^*F^{(2)} $ solves equation (\ref{Gznormal}) on $\mathcal M^{(1),rec} $.
Substituting these solutions in (\ref{extra}) and using (\ref{eq: Greens functions related})
 we see that
$G_{\mathcal O^{(i)}}(z; X,Y)$, $i=1,2$ satisfy
(\ref{eq: two greens functions}). Moreover, as 
the poles of $z\mapsto G_{\mathcal O}^{(i)}(z; X,Y)$ in $\C$,
that is the eigenvalues of $H_{\O^{(i)}}$,
 coincide for $i=1,2$, we see, 
using equations
(\ref{eq: jump of green's function}) and (\ref{eq: residues}), that
 BSP related to operators $H_{\O}^{(1)}$ on $\p \O^{(1)}$ 
and $H_{\O}^{(2)}$ on $\p \O^{(2)}$ are $\Phi^{rec}$-related. \qed
\medskip



\subsection{BSP for subdomains of $\mathcal M_{int}$ and recognition of singular points}
\label{subsec: Relation}
When $\Gamma\subset  \p \O^{(i)}$ and $s>0$, we denote the domain of influence by
$$
 \Omega^{(i)}_{\mathcal O}(\Gamma,s)=\{X\in \Omega^{(i)}_{\mathcal O};\
{\tilde d}^{(i)}(X,\Gamma)<s\}.
$$
where ${\tilde d}^{(i)}$ now is the distance in $\Omega^{(i)}_{\mathcal O}$.



\begin{theorem} \label{recognition}
Let $({\mathcal M}^{(1),rec},{\mathcal M}^{(2),rec},\Phi^{rec})$
be an admissible triple and 
${\mathcal O}^{(i)}=B^{(i)}(q_i,R) \subset {\mathcal M}^{(i),rec}$, $i=1,2,$
be a ball centered at $q_i$ and radius $R$ satisfying (\ref{eq: basic assumptions})
and (\ref{eq: basic assumptions2}). 
Denote 
\beq \label{5.22a}
s^{(i)}(q_i)=\min(d^{(i)}(q_i, {\mathcal M}^{(i)}_{sing}), (y_0-2)/4).
\eeq
Then $s^{(1)}(q_1)=s^{(2)}(q_2)$.
Using the notation $s=s^{(1)}(q_1)$, then, for 
\ba
   {\widetilde{\mathcal M}}^{(i),rec}=
 {\mathcal M}^{(i),rec} \cup \Omega^{(i)}_{\mathcal O}( \p \O^{(i)},s-R),
 \quad i=1,2,
  \ea
  there is a map
  $  {\widetilde \Phi}^{rec}:{\widetilde{\mathcal M}}^{(1),rec}
 \to
 {\widetilde{\mathcal M}}^{(2),rec}
$ 
which is an extension
of  $\Phi^{rec}$. Moreover, the triple 
$( {\widetilde{\mathcal M}}^{(1),rec}, 
{\widetilde{\mathcal M}}^{(2),rec},  {\widetilde \Phi}^{rec})$ is  admissible.
  \end{theorem}
\noindent {\ntekst Note that $B^{(i)}(q_i,s)=\Omega^{(i)}_{\mathcal O}( \p \O^{(i)},s-R)\cup B^{(i)}(q_i,R)$ and that,
$d^{(i)}(X, \p \O^{(i)})={\tilde d}^{(i)}(X, \p \O^{(i)})$ for  $X \in {\mathcal M}^{(i),rec} \setminus {\mathcal O}^{(i)}$.}

\noindent
{\bf Proof.} Assume opposite to the claim  that we would have
$s^{(1)}(q_1) > s^{(2)}(q_2)$.
Let 
\beq \label{5.21a}
a=s^{(1)}(q_1)-R,\quad b=s^{(2)}(q_2)-R, \quad 0<c<b<a.
\eeq
Then, by (\ref{decom of R map}), the BSP of the operator
$H_{\O}^{(i)}$ on $\p{\mathcal O}^{(i)}$ determines, on $\p \O^{(i)}$,  
the hyperbolic N-to-D map
$R^{(i), T}:=R^T_{\p \O^{(i)}}$ of the Riemannian surface 
$\Omega^{(i)}_{\mathcal O},$ $i=1,2$.
This and Lemma \ref{lem. determination of Green} 
yield that these maps satisfy
\beq \label{6.3.6}
(R^{(1), T}(h\circ \Phi^{rec}))(X)=
(R^{(2), T}(h))(\Phi^{rec}(X)),\quad X\in \p \O^{(1)},
\eeq
for all $h\in C^\infty_0(\p{\O}^{(2)} \times \R_+)$.

Let us deform the surfaces $\Omega^{(i)}_{{\mathcal O}}$ 
 replacing the metric with a smooth
metric in the $(b-c)/2-$neighborhood of the conic points and replacing the ends of the manifolds with compact surfaces.
We can do this by smoothly pinching
the first end-cylinder, $S^{r_1}\times (3/4 y_0+1/2, y_0)\subset \M_1$,  to a 
semisphere $S^2_+(r_1)$ and the parts of the other ends, $\M^{(i)}_j,\, j>1,$ 
which lie outside $\Omega^{(i)}_{\mathcal O}( \p \O^{(i)}, a)$,
 also to appropriate semispheres.
These give rise to two smooth
compact Riemannian surfaces ${\mathcal N}^{(i)}, i=1,2,$  with  
$\Gamma^{(i)}:=\p {\mathcal N}^{(i)}=\p\O^{(i)}$.
Then
the $c$-neighborhoods of $\p {\mathcal N}^{(i)}$ in ${\mathcal N}^{(i)}$, denoted by 
${\mathcal N}^{(i)}(\Gamma^{(i)}, c)$, are  isometric to 
 $\Omega^{(i)}_{\mathcal O}( \p \O^{(i)},c)$. By the finite velocity of the wave propagation,
 which is equal to one with respect to the underlying metric, the above isometry implies that 
 the N-to-D map $R^{(i),T}_{\p {\mathcal N}^{(i)}}$ on $\p {\mathcal N}^{(i)}$ 
 corresponding to manifold 
 ${\mathcal N}^{(i)}$ coincide with the N-to-D map on $\p \O^{(i)}$ 
 corresponding to manifold 
 $\Omega_\O^{(i)}$ for $T<2c$.
Together with (\ref{6.3.6}), this implies 
that the inverse of the N-to-D maps $R^{(i),T}_{\p {\mathcal N}^{(i)}}$, called the hyperbolic
D-to-N maps, satisfy the equation similar to (\ref{6.3.6}). 
By   \cite[Lemma 4.24 and p. 200]{KKL01}, the 
D-to-N maps with time $T<2c$,
determine uniquely the manifolds ${\mathcal N}^{(i)}(\p {\mathcal N}^{(i)},c)$, implying that
 there exists an isometry, 
 \ba
 {\widetilde \Phi}_c: {\mathcal N}^{(1)}(\p {\mathcal N}^{(1)},c) \to 
 {\mathcal N}^{(2)}(\p {\mathcal N}^{(2)},c).
 \ea

 Note that, if we identify $\p{\mathcal O}^{(1)}$ with $\p{\mathcal O}^{(2)}$, then
  the representation of this map in the boundary normal coordinates, see e.g.\ \cite{KKL01},
 is the identity map.
As ${\mathcal N}^{(i)}(\p {\mathcal N}^{(i)},c)$ is isometric to 
$ \Omega^{(i)}_{\mathcal O}(\p{\mathcal O}^{(i)},c)$ and above  $c<b$ is arbitrary, 
this implies that  there is an isometry
 \beq\label{eq: new isometry}
 {\widetilde \Phi}: \Omega^{(1)}_{\mathcal O}(\p{\mathcal O}^{(1)},b) \to
 \Omega^{(2)}_{\mathcal O}(\p{\mathcal O}^{(2)},b), 
 \eeq
{\ntekst By the conditions of Theorem, if $b' <b$ is so
   small that $\Omega_{\mathcal O}^{(i)}(\p {\mathcal O}^{(i)}), b') \subset
  {\mathcal M}^{(i),rec} $, then}
 \beq \label{5.25b} 
 {\widetilde \Phi}(X)=\Phi^{rec}(X), \quad X \in \Omega_{\mathcal O}^{(1)}(\p {\mathcal O}^{(1)}), b').
 \eeq

 As Green's functions  $G^{(i)}(z, X, Y)$, $i=1,2,$ satisfy
  relation (\ref{eq: Greens functions related})
 for 
$X, Y \in \Omega^{(i)}_{\mathcal O}( \p \O^{(i)},b')$,
we see, using the unique continuation in $X$ and $Y$ variables
as in the proof of Lemma \ref{lem. N1 admissible}, 
that 
\beq\label{eq: Greens functions related2}
& &G^{(2)}(z; {\widetilde \Phi}(X), {\widetilde \Phi} (Y))=G^{(1)}(z; X, Y),
\\ \nonumber
& &{for }
\ z\in \C\setminus \R, \quad
X,Y\in\Omega^{(1)}_{\mathcal O}(\p{\mathcal O}^{(1)},b).
\eeq 
Thus  $(\Omega^{(1)}_{\mathcal O}(\p{\mathcal O}^{(1)},b),
 \Omega^{(2)}_{\mathcal O}(\p{\mathcal O}^{(2)},b), {\widetilde \Phi})$
 is admissible. 
 {\ntekst Using (\ref{eq: new isometry}), (\ref{eq: Greens functions related2}), 
it follows from Remark \ref{Remark:3} that $\Phi^{rec}$ can be extended by ${\widetilde \Phi}$
as ${\widetilde \Phi}^{rec}$,
\beq \label{5.23a}
& &{\widetilde \Phi}^{rec}: {\widetilde \M}^{(1), rec} \rightarrow
{\widetilde \M}^{(2), rec};
\\\nonumber
& & {\widetilde \M}^{(i), rec}= \M^{(i), rec} 
\cup \Omega^{(i)}_{\mathcal O}(\p {\mathcal O}^{(i)},b)
\eeq

}

Recall that, by our assumption, $a>b$.
{\nntext Due to  (\ref{eq: basic assumptions}), (\ref{5.22a}), this implies that
${\mathcal M}^{(2)}_{sing}\cap \p\Omega^{(2)}_{\mathcal O}(\p{\mathcal O}^{(2)},b)\not =\emptyset$.
Next we show that this is not possible.
 

For $Y\in  \p \O^{(i)}$ 
 we define the  boundary-cut-locus distance 
\ba
\tau_{{\O}^{(i)}}(Y)=\inf\{t>0;\ 
\gamma_{Y,\nu}^{(i)}(t)\in  {\mathcal M}^{(i)}_{sing}\hbox{ or }
 {\tilde d}^{(i)}(\gamma_{Y,\nu}^{(i)}(t), \p \O^{(i)})<t\},
\ea
where 
$\nu \in T_Y {\mathcal M}^{(i)}$ is the exterior unit normal vector to 
$ \p \O^{(i)}$ and $\gamma^{(i)}_{Y,\nu}(t)$ is  the geodesic on ${\mathcal M}^{(i)}_{int}$.

As the mapping
(\ref{eq: new isometry}) is an isometry between
$ \Omega^{(1)}_{\mathcal O}(\p{\mathcal O}^{(1)},b)$
and $ \Omega^{(2)}_{\mathcal O}(\p{\mathcal O}^{(2)},b)$,
we see that, for all $Y\in \p{\mathcal O}^{(1)}$,
\ba
\min (\tau_{{\O}^{(1)}}(Y),b)=\min (\tau_{{\O}^{(2)}}(\Phi^{rec}(Y)),b).
\ea

Next, any point $p^{(i)}\in \Omega^{(i)}_{\mathcal O}(\p{\mathcal O}^{(i)},b)$ 
 can be written in the form
 $\gamma^{(i)}_{Y, \nu}(t)$ where $Y\in \p \O^{(i)}$ and 
 $t  \leq 
 \min (\tau_{{\O}^{(i)}}(Y),b)$.
  Moreover, if 
 \beq \label{5.24a}
 p^{(i)} \in \p\Omega^{(i)}_{\mathcal O}(\p{\mathcal O}^{(i)},b) \setminus \p {\mathcal O}^{(i)},
 \eeq 
 then $d^{(i)}(p^{(i)}, \p {\mathcal O}^{(i)})=b$ and
 there is 
 \beq \label{5.24b}
 Y^{(i)} \in \p{\mathcal O}^{(i)}\quad \hbox{such\,that}\quad b \leq \tau_{{\O}^{(i)}}(Y^{(i)}),
 \quad
 p^{(i)}=\gamma^{(i)}_{Y^{(i)}, \nu}(b). 
 \eeq
 
  
Let 
\beq \label{5.??}
p^{(2)} \in {\mathcal M}^{(2)}_{sing}\cap \p\Omega^{(2)}_{\mathcal O}(\p{\mathcal O}^{(2)},b).
\eeq
By  definition (\ref{5.22a}), (\ref{5.21a}), $p^{(2)}$
satisfies (\ref{5.24a}), (\ref{5.24b}).
By the above, there is 
a point $Y^{(2)}\in  \p \O^{(2)}$ such that $p^{(2)}= \gamma^{(2)}_{Y^{(2)}, \nu}(b)$.
Let $Y^{(1)}=(\Phi^{rec})^{-1}(Y^{(2)})$ and consider $p^{(1)}=\gamma^{(1)}_{Y^{(1)}, \nu}(b)$.
Since $\Omega^{(1)}_{\mathcal O}(\p{\mathcal O}^{(1)},b)$  and 
$\Omega^{(2)}_{\mathcal O}(\p{\mathcal O}^{(2)},b)$ are isometric, $p^{(1)}$
satisfies (\ref{5.24a}), (\ref{5.24b}). Moreover, since $a<b$,
$p^{(1)} \notin \M^{(1)}_{sing}.$

Let 
\ba
p^{(i)}_{\e}:= \gamma_{Y^{(i)}, \nu}(b-2\e),\quad \e >0,\ i=1,2.
\ea

For $\e < b/8$, denote by ${\widetilde O}^{(i)}_\e= B^{(i)}(p^{(i)}_{\e}, \e)$  the 
metric ball in  $\Omega^{(i)}_{\mathcal O}$ of radius $\e$.
By using (\ref{5.23a}) and choosing
$\e>0$ to be  small,  ${\widetilde \O}^{(i)}_\e$ satisfy the conditions of Lemma
\ref{lem. determination of Green} with ${\widetilde \M}^{(i)}$ instead of $\M^{(i)}$,
${\widetilde \Phi}^{rec}$ instead of $\Phi^{rec}$
and ${\widetilde \O}^{(i)}_\e$ instead of $\O^{(i)}$.

Then, Lemma \ref{lem. determination of Green} implies that
\beq \label{7.3.2}\quad \quad
\hbox{BSP \, for}\,  H_\e^{(i)} \,\, \hbox{on} \,\,{\widetilde \O}_\e^{(i)},\,\, i=1,2,\,\,
\hbox{are}\,\, 
{\widetilde \Phi}^{reg}\hbox{-related}.
\eeq
Here $H_\e^{(i)}=H^{(i)}_{{\widetilde \O}_\e}$, i.e. is the Laplace operator
associated with ${\widetilde \Omega}^{(i)}_\e=\M^{(i)}_{int} \setminus {\widetilde \O}^{(i)}_\e$. 
Equation (\ref{7.3.2}) together with Lemma \ref{lem: volumes} imply that
\ba
S^{(1)}({\widetilde \Omega}_\e^{(1)}(\p {\widetilde \O}_\e^{(1)},r-\e))=
S^{(2)}({\widetilde \Omega}_\e^{(2)}(\p {\widetilde \O}_\e^{(2)},r-\e))
\ea
when $r>0$. }

Since ${\widetilde \Phi}^{rec}$ is an isometry, we also have
\ba
S^{(1)}({\widetilde \O}_\e^{(1)})=S^{(2)}({\widetilde \O}_\e^{(2)}).
\ea
On the other hand, when $\e>0$ is small enough, 
\ba
S^{(i)}(B^{(i)}(p^{(i)}_\e, r))=S^{(i)}({\widetilde \O}_\e^{(i)})+
S^{(i)}({\widetilde \Omega}_\e^{(i)}(\p {\widetilde \O}_\e^{(i)},r-\e)).
\ea
Therefore, the above two equations imply that
\beq \label{7.3.3}
S^{(1)}(B^{(1)}(p^{(1)}_\e, r))=S^{(2)}(B^{(2)}(p^{(2)}_\e, r)).
\eeq
Next, we observe that as
$d^{(i)}(p^{(i)}_{\e}, p^{(i)}) \leq 2\e,$ we have
$$
B^{(i)}(p^{(i)},r-2\e) \subset B^{(i)}(p^{(i)}_{\e},r) \subset 
B^{(i)}(p^{(i)},r+2\e),\quad \hbox{for }r>2\e.
$$
Thus,
by the continuity of the area,
\ba
S^{(i)}(B^{(i)}(p^{(i)},r))
 = \lim_{\e \to 0} S^{(i)}(B^{(i)}(p^{(i)}_{\e},r)).
\ea
Together with (\ref{7.3.3}), this implies that, for $r>0$,
\beq \label{E2}
S^{(1)}(B^{(1)}(p^{(1)},r))=S^{(2)}(B^{(2)}(p^{(2)},r)).
\eeq
Let us now consider the polar coordinates of $\mathcal M^{(i)}$ 
near $p^{(i)}$ where we note that, due to 
$d^{(i)}({\mathcal O}^{(i)}, \Gamma_0) >(y_0-2)/2$, we have $p^{(i)} \notin \Gamma_0$. 
In these coordinates,
$$
(ds)^2 = (dr)^2 + C^{(i)}r^2\left(1 + h^{(i)}(r,\theta) \right)(d\theta)^2,
$$
cf.\ (\ref{S1metricaroundsingular}). It then follows from (\ref{E2}), that
$$
 C^{(i)}=C^{(i)}(p^{(i)})= 
 \lim_{r \to 0} \frac{1}{\pi r^2} \left[S^{(i)}(B^{(i)}(p^{(i)},r) )    \right] 
 $$
satisfy 
\beq \label{7.3.5}
C^{(1)}(p^{(1)})=C^{(2)}(p^{(2)}).
\eeq
Note that, 
if $p^{(i)}\in {\mathcal M}_{sing}^{(i)}$ we have $C^{(i)}\neq 1$
and if
  $p^{(i)}\in {\mathcal M}_{reg}^{(i)}$ then $C^{(i)}=1$.
  As we assume that $a=s^{(1)}(q_1)-R>b=s^{(2)}(q_2)-R$, we have 
   $p^{(1)}\in {\mathcal M}_{reg}^{(1)}$ 
   and thus $C^{(1)}=1$. Hence, we also have   $C^{(2)}=1$, 
 and thus $p^{(2)}\in  {\mathcal M}_{reg}^{(2)}$, contradicting (\ref{5.??}).
%
This implies that $a \leq b$ which is in contradiction with our assumption
that we would have
$s^{(1)}(q_1) > s^{(2)}(q_2)$. This shows that we must have
\ba
s^{(1)}(q_1)=s^{(2)}(q_2).
\ea

This equation together with (\ref{5.23a}) prove the theorem.
\qed

\medskip
%

Let $\mathcal A$ be the collection of admissible triples
$({\mathcal W}^{(1)},\, {\mathcal W}^{(2)}, \Phi  )$ such that 
$N \subset {\mathcal W}^{(1)},\, i=1, 2$.
We define a partial order on $\mathcal A$ by setting
$({\mathcal W}^{(1)},\, {\mathcal W}^{(2)}, \Phi  )\leq
(\widetilde {\mathcal W}^{(1)},\, \widetilde{\mathcal W}^{(2)},\widetilde \Phi  )$
if ${\mathcal W}^{(1)}\subset \widetilde {\mathcal W}^{(1)}$ 
{\ntekst and
$
\Phi= {\widetilde \Phi}|_{{\mathcal W}^{(1)}}.
$

Note that, by Remark \ref{Remark:3}, if 
$({\mathcal W}^{(1)},\, {\mathcal W}^{(2)}, \Phi  )$ and 
$(\widetilde {\mathcal W}^{(1)},\, \widetilde{\mathcal W}^{(2)},\widetilde \Phi  )$
are admissible triples, then $({\mathcal W}^{(1), en},\, {\mathcal W}^{(2), en}, \Phi^{en}  )$,
where
\ba
& &{\mathcal W}^{(i), en}={\mathcal W}^{(i)} \cup {\widetilde {\mathcal W}}^{(i)},
\\ \nonumber
& & \Phi^{en}|_{{\mathcal W}^{(1)}}=\Phi, \quad
\Phi^{en}|_{{\widetilde {\mathcal W}}^{(1)}}={\widetilde \Phi},
\ea
is also an admissible triple. }
Therefore,  by Zorn's lemma, there exists a maximal element 
$({\mathcal W}_m^{(1)},\, {\mathcal W}_m^{(2)}, \Phi_m  ) \in \mathcal A$.


\begin{lemma} \label{maximal}
The maximal element $({\mathcal W}_m^{(1)},\, {\mathcal W}_m^{(2)}, \Phi_m  ) $ of $\mathcal A$ satisfies
 \beq \label{10.3.6}
 {\mathcal W}_m^{(1)}= \M^{(1)}_{reg}.
 \eeq 
 \end{lemma}
 \noindent
{\bf Proof.} If the claim is not true, there exists
    $X^{(1)}_0 \in\ M ^{(1)}_{reg} \cap \p{\mathcal W}_m^{(1)}$.
    Let $\mu([0,1])$ be a smooth path from 
    $\mu(0)=Z =(x, y),\, x \in \Gamma_0,\, y=2/3+y_0/3$ to $\mu(1)=X^{(1)}_0$, such that 
    \ba
    \mu([0,1)) \subset \M^{(1)}_{reg},\quad \mu \cap \left(\Gamma_0 \times (\frac{2+y_0}{2}, y_0) \right) =\emptyset.
    \ea
    Then $d_0=d^{(1)}(\mu, \M^{(1)}_{sing})>0$. Let 
     $c= \min(\frac{y_0-2}{4}, \, \frac{d_0}{2})$.  
    We can cover $\mu([0, 1])$ by a finite number of balls $B^{(1)}_j=B^{(1)}(X^{(1)}_j, c/2) \subset 
    \M ^{(1)}_{reg}$ 
    so that  
    \beq \label{10.3.5}
   { \overline B}_j^{(1)} \subset {\mathcal W}_m^{(1)},\,B_j^{(1)} \cap \Gamma_0=\emptyset, \,\,
    X^{(1)}_{j+1} \in B_j^{(1)},
   \eeq 
   where we order them so that $X^{(1)}_0 \in B^{(1)}_1$.
   Let  ${\mathcal O}^{(1)}_1=B^{(1)}(X^{(1)}_1, R)$  be a small ball such 
   that $0<R<c/2$ satisfies  (\ref{eq: basic assumptions}), (\ref{eq: basic assumptions2}),
   and 
   ${\mathcal O}^{(1)}_1   \subset {\mathcal W}_m^{(1)}$. 
   As $d^{(1)}( X^{(1)}_1,{\mathcal M}^{(1)}_{sing})>\frac {d_0}2$,
    Theorem \ref{recognition} yields that
 {\ntekst we can extend 
 the admissible triple $({\mathcal W}_m^{(1)},\, {\mathcal W}_m^{(2)}, \Phi_m  ) $ onto
 \ba
 {\widetilde {\mathcal W}}^{(i)}={\mathcal W}_m^{(i)} \cup B^{(i)}(X^{(i)}_1, c),
 \quad X^{(2)}_1= \Phi_m(X^{(1)}_1).
 \ea}
 As $ X^{(1)}_0\in B( X^{(1)}_1, c)$, this contradicts 
the fact that $({\mathcal W}_m^{(1)},\, {\mathcal W}_m^{(2)}, \Phi_m  ) $ is a maximal element
 of $\mathcal A$, which completes the proof of (\ref{10.3.6}).
    \qed\medskip

   Lemma \ref{maximal} proves that there is a diffeomorphism
   \ba
\Phi_m: {\mathcal M} ^{(1)}_{reg} \to {\mathcal W}_m^{(2)}, \quad 
{\mathcal W}_m^{(2)}= \Phi_m({\mathcal M} ^{(1)}_{reg})
\subset  {\mathcal M} ^{(2)}_{reg},
\ea
which is  a Riemannian isometry.
   Changing the role of indexes 1 and 2, we see that  there is also a  diffeomorphism
\ba
{\widetilde \Phi}_m: {\mathcal M} ^{(2)}_{reg} \to {\widetilde {\mathcal W}}_m^{(1)}, 
\quad {\widetilde {\mathcal W}}_m^{(1)}
\subset  {\mathcal M} ^{(1)}_{reg}
\ea
which is a Riemannian isometry. Moreover,
using Lemma \ref{lem. N1 admissible} we see that ${\widetilde \Phi}_m$
and $\Phi_m$ coincide with the identity map on $\Gamma_0$.

Using (\ref{eq: Greens functions related}) we see that for all  $z\in \C\setminus \R$,
$X\in {\mathcal M} ^{(2)}_{reg} $ and $Y\in \Gamma_0$.
\ba
G^{(1)}(z; \Phi_m({\widetilde \Phi}_m(X)), Y)=G^{(2)}(z; X, Y).
\ea
 By Lemma \ref{identification}, this implies that $\Phi_m({\widetilde \Phi}_m(X))=X$,
 that is, $\Phi_m\circ{\widetilde \Phi}_m=I$ on ${\mathcal M} ^{(1)}_{reg}$.
 Similarly, we see that   $\widetilde \Phi_m\circ{ \Phi}_m=I$
 on ${\mathcal M} ^{(2)}_{reg}$ and hence
 \ba
 {\mathcal W}_m^{(2)}={\mathcal M} ^{(2)}_{reg},\,\, {\mathcal W}_m^{(1)}={\mathcal M} ^{(1)}_{reg},
 \quad \hbox{and}\,\, {\widetilde \Phi}_m= \Phi_m^{-1}.
\ea 
Summarizing, we have shown that there is a diffeomorphism
\ba
\Phi_m:({\mathcal M}_{reg}^{(1)}, g^{(1)}) \to ({\mathcal M}_{reg}^{(2)}, g^{(2)}),
\ea 
which is a Riemannian isometry.

{\nntext 
Skipping again the superscript $^{(i)}$, we show next that
\beq\label{distances are same}
d(X, Y)= d_{reg}(X, Y),\quad \hbox{for any $X, Y \in {\mathcal M}_{reg}$},
\eeq
where $d_{reg}$ is the distance on $({\mathcal M}_{reg}, g)$ defined
as the infimum of the length of  rectifiable paths connecting $X$ to $Y$.
As ${\mathcal M}_{reg}\subset {\mathcal M}$, we have
$d(X, Y)\leq d_{reg}(X, Y)$. On the other hand, let
$X,Y\in {\mathcal M}_{reg}$ and consider a  rectifiable path $\mu:[0,\ell]\to {\mathcal M}$  
from $X$ to $Y$,
parametrized by the arc-length.
As we consider infimum of the length of paths, we can assume that $\mu$ is 
one-to-one. As ${\mathcal M}_{sing}$ is discrete, $\mu$ can intersect it only
finite many times. If $p=\mu(t_0)\in {\mathcal M}_{sing}$, let us consider
the coordinates $X:U\to [0,\e_p)\times [0,2\pi]$ near $p$ defined in (A-2). Let  $\e>0$ be small enough
and $t_-,t_+\in (0,\ell)$, $t_-<t_0<t_+$ be such that $X(\mu(t_\pm))=(\e,\theta_\pm)$. 
If we then modify the path $\mu$
by replacing $\mu([t_-,t_+])$ by a segment on the circle, that is, the path $X^{-1}(\e,J)$ where $J\subset [0,\pi]$
in the interval connecting $\theta_-$ to $\theta_+$,  the length of $\mu$
is increased by $O(\e)$. By choosing $\e$ small enough and modifying
the path $\mu$ in the above way in all points where $\mu$ intersects ${\mathcal M}_{sing}$,
we see that near $\mu$ there is a path in ${\mathcal M}_{reg}$ which length
is arbitrarily close to the length of $\mu$. This shows that $d(X, Y)\geq d_{reg}(X, Y)$
proving (\ref{distances are same}).

The identity  (\ref{distances are same}) implies that $
({\mathcal M}_{int}, d)$, considered as a metric space, is isometric to
 the completion of the metric space 
$({\mathcal M}_{reg},\, d_{reg}).$}
Thus, we can uniquely extend $\Phi_m$ to a metric isometry
\begin{equation} \label{E3B}
 \Phi:({\mathcal M}^{(1)}_{int}, d^{(1)}) \to ({\mathcal M}^{(2)}_{int}, d^{(2)}),
\end{equation} 
Again, taking into account that the number of singular points is finite,
we see that $ \Phi$ maps singular points to singular points.
Let us numerate the singular points on ${\mathcal M}^{(1)}$ and ${\mathcal M}^{(2)}$
as $p_l^{(1)}\in {\mathcal M}^{(1)}$,  $l=1,2,\dots,L$ and 
 $p_l^{(2)}\in {\mathcal M}^{(2)}$,  $l=1,2,\dots,L$ so that $p_l^{(2)}=  \Phi_m(p_l^{(1)})$.
%

The map $\Phi$ defined above  satisfies conditions (1)--(3) of Theorem \ref{thm: main result}.
We  prove (4) we use the following Lemma:


 \begin{lemma} \label{singular}
Let 
\beq \label{11.3.3}
\Phi: {\mathcal M}^{(1)} \to {\mathcal M}^{(2)}
\eeq 
satisfy the conditions (1)--(3) of Theorem \ref{thm: main result}. Then $\Phi$ satisfies the
condition (4).
 \end{lemma}
 \noindent
{\bf Proof.}
Let $p^{(1)}_l \in \M^{(1)}_{sing},\,p_l^{(2)}=\Phi(p_l^{(1)}) \in \M^{(2)}_{sing}$
and $\e_0>0$ be so small that polar coordinates  (A-2)
centered at $p^{(i)}_l$ are well defined in the ball $B^{(i)}(p^{(i)}_l,\e_0)$
for $i=1,2$. We denote these coordinates by $\psi^{(i)}:B^{(i)}(p_l^{(i)},\e_0)
\to  [0,\e_0)\times [0,2\pi)$. Below, we skip the subindex $l$.

First, using a point $q_1\in \M^{(1)}_{reg}$  such that $p^{(1)}$
is the unique closest singular point of $\M^{(1)}$ to $q_1$,
the proof of Theorem \ref{recognition}, see (\ref{7.3.5}), shows that
\beq \label{11.3.6}
C_1=C^{(1)}(p^{(1)}) =C^{(2)}(p^{(2)})=C_2:=C.
\eeq

Let us consider a distance minimizing curve in $B^{(i)}(p^{(i)},\e_0)$ emanating  
from the point $p^{(i)}$. We call such curve a radial geodesics 
and denote it by $\gamma^{(i)}(s)$ where $s$ is the arclength from $p^{(i)}$.

By (\ref{S1metricaroundsingular})
the radial geodesic $\gamma^{(i)}=\gamma^{(i)}([0,\e_0))$ is given 
in normal coordinates by 
$\psi^{(i)}(\gamma^{(i)})=\{(\theta,r);\ \theta=\a_0,\ 0\leq r<\e_0\}$ where $\a_0\in [0,2\pi)$
is a parameter associated to $\gamma^{(i)}$, and we denote below 
 $\a^{(i)}(\gamma^{(i)})=\a_0$.

Since $\Phi$ is an isometry, it maps
any radial geodesic $\gamma^{(1)}$  emanating  from $p^{(1)}$
to some radial geodesic  $\gamma^{(2)}$ emanating  from $p^{(2)}$.  
When $\gamma^{(i)}_0$ is the geodesic
satisfying $\a^{(i)}(\gamma_0^{(i)})=0$,
the  parameter $\a^{(i)}(\gamma^{(i)})$ associated to $\gamma^{(i)}(s)$ 
satisfies
%
\ba
\a^{(i)}(\gamma^{(i)})= 
\lim_{\e \to 0}\frac{\ell^{(i)}_\e(\gamma^{(i)}(\e),\,\gamma^{(i)}_0(\e))}{C_i \e}, 
\ea
where $\ell^{(i)}_\e(\gamma^{(i)}(\e),\,\gamma^{(i)}_0(\e))$ is the arc length of the counter-clockwise  oriented path 
connecting $\gamma^{(i)}_0(\e)$ and $\gamma^{(i)}(\e)$ along the circle 
$S^{(i)}_\e=\{X\in \mathcal M^{(i)};\,d^{(i)}(X, p^{(i)})=\e\}$. 
Let $\beta=a^{(1)}(\gamma^{(1)})-\a(\gamma^{(2)})\in (-2\pi, 2\pi)$.
As $\Phi$ is an isometry,
\beq \label{5.36}
h_l^{(1)}(r, \theta)=h_l^{(2)}(r, {\widehat{\theta+\beta}}),
\eeq 
where, for $\theta \in \R,\, {\widehat \theta} \in [0, 2\pi)$ satisfies 
$\theta -{\widehat \theta} \in 2 \pi \Z$.

This completes the proof of Lemma \ref{singular} and Theorem \ref{thm: main result}.\qed


\section{Orbifold isomorphism for $\Gamma\backslash{\bf H}^2$}
We shall prove Theorem 1.3. Let $\mathcal M^{(i)} = \Gamma_i\backslash{\bf H}^2$ 
and 
$\mathcal M^{(i)}_{sing}$ the set of elliptic singular points in $\mathcal M^{(i)}$. We 
have already constructed a hyperbolic isometry 
$\Phi : \mathcal M^{(1)}\setminus{\mathcal M}^{(1)}_{sing} \to \mathcal M^{(2)}\setminus\mathcal M^{(2)}_{sing}$ in \S 5. Since the hyperbolic metric is conformal to the Euclidean metric, 
$\Phi$ is conformal. As $\mathcal M^{(i)}$ is orientable, we can assume 
$\Phi : \mathcal M^{(1)}\setminus{\mathcal M}^{(1)}_{sing} \to 
\mathcal M^{(2)}\setminus\mathcal M^{(2)}_{sing}$ to be analytic. Take 
$p^{(1)} \in \mathcal M^{(1)}_{sing}$ and a small disc $B^{(1)}(p, {\epsilon})$ centered 
at $p^{(1)}$. 
Since $\Phi$ maps $p^{(1)}$ to $p^{(2)}:=\Phi(p^{(1)}) \in \mathcal M^{(2)}_{sing}$, 
$p^{(1)}$ is a removable 
singularity for $\Phi$. Hence $\Phi$ is analytic also at $p^{(1)}$.
Let $({\widetilde B}^{(i)}(P^{(i)}, {\epsilon}),{\widetilde g}^{(i)})$ be the
uniformizing covers of  $({B}^{(i)}(p^{(i)}, {\epsilon}),{g}^{(i)})$, $i=1,2$.
Then $\Phi$ can be lifted to an analytic map between the coverings except for the center  
$$
\widetilde \Phi : {\widetilde B}^{(1)}(P^{(1)},{\epsilon})\setminus\{P^{(1)}\} \to 
{\widetilde B}^{(2)}(P^{(2)}, {\epsilon}\setminus\{P^{(2)}\}.
$$

 This implies that $P^{(1)}$ is a removable singularity of 
 $\widetilde \Phi$, hence $\widetilde \Phi$ is analytic on 
 ${\widetilde B}^{(1)}(P^{(1)},\epsilon)$.
  
 {\ntekst It follows from (\ref{5.36}) that
  $p^{(1)}$ is a singular point of the orbifold if and only if
  $ p^{(2)}$ is a singular point. Moreover, the map 
 \ba
 r \mapsto r,\quad \frac{\theta}{n} \mapsto \frac{{\widehat{\theta+\beta}}}{n},
 \,\, n=C^{-1/2},
 \ea
 extends to the isometry between ${\widetilde B}^{(i)}(P^{(i)}, {\epsilon}),\, i=1,2$.
 }
 This completes the proof of Theorem 1.3. \qed

\medskip
By the suitable linear transformation $\gamma \in SL(2,{\bf R})$, $\mathcal M_{\nu_1}^{(1)}$ is mapped to $\mathcal M_{\nu_2}^{(2)}$ conformally. Identifying them, we see that $\Phi$, constructed above, is the identity on $\mathcal M_{\nu_1}^{(1)}$, hence is equal to the identity on all of $\mathcal M^{(1)}$. This implies that $\gamma \Gamma_1 \gamma^{-1}\backslash{\bf H}^2 = \Gamma_2\backslash{\bf H}^2$, hence $\gamma\Gamma_1\gamma^{-1} = \Gamma_2$. Therefore, the generalized S-matrix determines the conjugate class of geometrically finite Fuchsian groups.

\begin{remark} \label{rem:boundary}
The technique used in this paper can be easily extended to consider the case when
$\p \M \neq \emptyset$. In this case we should require, in addition to  {\bf (A-1)},
{\bf (A-2)}, that each end is diffeomorphic to either a cylinder or a strip 
$(0, \ell_i) \times (1, \infty)$ with   the metric satisfying {\bf (A-3)},
{\bf (A-4)} where, in the case of a strip, $0 \leq x \leq \ell_i$.
\end{remark}

 \subsection*{Acknowledgments}
The research of Y.K. is partially
supported by EPSRC Grant EP/F034016/1 and Mathematical Sciences
Research Institute (MSRI). The research of M.L.
was financially supported by the Academy of Finland Center of Excellence programme 213476
and MSRI. The research of H.I. was supported by Gakusin Grant in Aid No. 21340028.

\end{document}